\documentclass[11pt]{elsarticle}
\usepackage{epsf,exscale,times}
\usepackage{amsmath,amsfonts,amssymb,bm,graphicx}
\usepackage[usenames]{color}

\usepackage{eurosym}

\usepackage{listings}
\usepackage{subfigure}
\usepackage{multicol}
\usepackage{multirow}
\usepackage{color}

\journal{Mathematics and Computers in Simulation}

\usepackage[colorlinks=true, linkcolor=blue, filecolor=blue, citecolor=blue, urlcolor=blue]{hyperref}

\usepackage{color}
\definecolor{gray97}{gray}{.97}
\definecolor{gray75}{gray}{.75}
\definecolor{gray45}{gray}{.45}
\definecolor{gray35}{gray}{.35}

\usepackage{color}
\definecolor{gray97}{gray}{.97}
\definecolor{gray75}{gray}{.75}
\definecolor{gray45}{gray}{.45}
\definecolor{gray35}{gray}{.35}

\usepackage{listings}
\lstnewenvironment{listing}[1][]
   {\lstset{#1}\pagebreak[0]}{\pagebreak[0]}

\lstdefinestyle{consola}
   {basicstyle=\scriptsize\bf\ttfamily,
    backgroundcolor=\color{gray75},
   }

\lstdefinestyle{C}
   {language=C,
   }

\begin{document}
\title{\Large \bf Static and dynamic SABR stochastic volatility models: calibration and option pricing using GPUs  \tnoteref{t1}}
\tnotetext[t1]{Partially financed by MCINN (Project MTM2010--21135--C02-01) and Ayuda CN2011/004 partially funded with FEDER funds). Also the authors would like to thank M. Men\'endez, from Banesto Bank, for providing the market data.}

\author[uam]{J.L. Fern\'andez}
\ead{joseluis.fernandez@uam.es}

\author[udc]{A.M. Ferreiro}
\ead{aferreiro@udc.es}

\author[udc]{J.A. Garc\'{\i}a}
\ead{jagrodriguez@udc.es}

\author[udc]{A. Leitao}
\ead{alvaro.leitao@udc.es}

\author[udc]{J.G. L\'opez-Salas}
\ead{jose.lsalas@udc.es}

\author[udc]{C. V\'azquez}
\ead{carlosv@udc.es}

\address[uam] {Department of Mathematics, Universidad Aut\'onoma de Madrid, \\ Francisco Tom\'as y Valiente, 7, 28049 -- Madrid (Spain)}

\address[udc]  {Departament of Mathematics, University of A Coru\~na, \\ Campus Elvi\~na s/n, 15071 -- A Coru\~na (Spain)}

\date{}

\begin{abstract}
For the calibration of the parameters in static and dynamic SABR stochastic volatility models, we propose the application of the GPU technology to the Simulated Annealing global optimization algorithm  and to the Monte Carlo simulation.  This calibration has been performed for EURO STOXX 50 index and EUR/USD exchange rate with an asymptotic formula for volatility or Monte Carlo simulation. Moreover, in the dynamic model we propose an original more general expression for the functional parameters, specially well suited for the EUR/USD exchange rate case. Numerical results illustrate the expected behavior of both SABR models and the accuracy of the calibration. In terms of computational time, when the asymptotic formula for volatility is used the speedup with respect to CPU computation is around $200$ with one GPU. Furthermore, GPU technology allows the use of Monte Carlo simulation for calibration purposes, the computational time with CPU being prohibitive.
\end{abstract}
\begin{keyword}
Parallel Simulated Annealing, SABR volatility model, Calibration, GPUs, CUDA.
\end{keyword}
\maketitle

\section{Introduction}
Mathematical models have become of great importance in order to price financial derivatives on different underlying assets. However, in most cases there is no explicit solution to the governing equations, so that accurate robust fast numerical methods are required. Furthermore, financial models usually depend on many parameters that need to be calibrated to market data (market data assimilation). As in practice the model results are required almost at real time, the speed of numerical computations becomes critical and this calibration process must be performed as fast as possible.

In the classical Black-Scholes model \cite{Black-73}, the underlying asset follows a lognormal process with constant volatility.  However, in real markets the volatilities are not constant and they can vary for each maturity and strike (volatility surface). In order to overcome this problem, different local and stochastic volatility models have been introduced (see \cite{Dupire-94,Dupire-97, Hagan, Heston, Hull-White}, for example). Here we consider the SABR model first proposed in \cite{Hagan}, where a first order approximation formula for the implied volatility of European plain-vanilla options with short maturities is obtained. In \cite{obloj} this formula is improved. In \cite{Paulot2009} a general method to compute a Taylor expansion of the implied volatility is described. In particular, a second order asymptotic SABR volatility formula is also computed. Next, in \cite{TakashakiTakeharaToda2011} a fifth order asymptotic expansion is proposed for $\lambda-$SABR model; thus providing an extension with a mean-reversion term in \cite{Henry-Labordere2008}.

The existence of closed-form formula simplifies the calibration of the parameters to fit market data. However, when considering constant parameters (static SABR model), the volatility surface of a set of market data for several maturities cannot be suitably fitted. In \cite{Hagan,west}, the calibration of the static SABR model to fit a single volatility smile is analyzed. Among the different techniques to deal with a set of different maturities (see \cite{book:Gatheral-06,homescu2011}, for example), in \cite{Hagan} a SABR model with time dependent parameters (dynamic SABR) is introduced and in \cite{Osajima-2007} an asymptotic expression for the implied volatility is obtained. Also by means of piecewise constant parameters, in \cite{GlassermanWu2010} the static SABR model is extended. In \cite{art:Larson2010} a second order approximation to call options prices and implied volatilities is proposed and a closed form approximation of the option price extending dynamically the original SABR model is gained. In \cite{KaisajunttiKennedy2011} the SABR model with a time-dependent volatility function and a mean reverting volatility process is obtained. In \cite{ChenGrzelakOsterlee2011} a hybrid SABR$-$Hull-White model for long-maturity equity derivatives is considered.

However, time dependent parameters highly increase the computational cost and it is not always possible to compute an analytical approximation for the implied volatility or the prices (or the expression results to be very complex). In this case, we can use numerical methods (for example, Monte Carlo) in the calibration process. In order to calibrate a model, an efficient, robust and fast optimization algorithm has to be chosen, either a local optimization algorithm (such as Nelder-Mead or Levenberg-Marquardt ones) or a global optimization one (such as Simulated Annealing, genetic or Differential Evolution based ones). Although local optimization algorithms are efficient, if the calibration function presents several local minima they can stuck in any of these ones (note that the volatility surface can be uneven for some markets). On the other hand, global optimization algorithms are more robust, although they involve greater computational cost and are much slower. As financial instruments analysis should be carried out almost in real-time, we use the efficient implementation of the Simulated Annealing (SA) algorithm in Graphics Processing Units (GPUs) proposed in \cite{art:SA-GPU}. Particularly, we consider Nvidia Fermi GPUs and the API for its programming, named CUDA (Compute Unified Device Architecture), see \cite{book-nvidia-2011,NVIDIA-2011}. CUDA consists of some drivers for the graphic card, a compiler and a language that is basically a set of extensions for the C/C++ language. This framework allows to control the GPU (memory transfer  operations, work assignment to the processors and threads synchronization).

Once the parameters have been calibrated, the model can be used to price exotic options. There are different techniques for pricing, such as Monte Carlo simulation, finite difference methods, binomial trees or integration-based methods (Fourier transform, for example). Among them, Monte Carlo simulation is a flexible and powerful tool that allows to price complex options (see \cite{book:Glasserman2003,book:Jackel2002}). From the computational viewpoint, it results to be very expensive mainly due to its slow convergence. Once again this is a handicap, particularly in pricing and risk analysis in the financial sector. However, as illustrated in the present paper, Monte Carlo simulation is suitable for pricing options on GPUs \cite{art:LeeYauGilesDoucetHolmes_10}.

In this work, we have parallelized the Monte Carlo method on GPUs for the static and dynamic SABR models. In the literature the implementation of efficient Monte Carlo methods for pricing options has been analyzed. In \cite{Joshi2009} Asian options pricing with Black-Scholes model by a quasi-Monte Carlo method is considered, thus getting a speedup factor up to $150$. In \cite{art:Bernemann-Schereyer-11} a Monte Carlo GPU implementation for the multidimensional Heston and hybrid Heston$-$Hull-White models (using a hybrid Taus$-$Mersenne-Twister random number generator) is presented, achieving speedup factors around $50$ for Heston model, and from $4$ to $25$ for the Heston$-$Hull-White model. In \cite{TianZhuKlebanerHamza2010}, European and American option pricing methods on GPUs under the static SABR model are presented. For the European case a speedup of around $100$ is obtained, while for American ones it is around $10$. In \cite{TianZhuKlebanerHamza2011} the pricing of barrier and American options using a parallel version of least squares Monte Carlo algorithm is carried out, following the techniques presented in \cite{TianZhuKlebanerHamza2010}. They obtain speedups up to $134$ in the case of barrier options and around $22$ in American ones. In our present  work, the random number generation is performed ``on the fly'' by using the Nvidia CURAND library (see \cite{NVIDIA-2011}, for details). This is an important difference with previous works, where random numbers are previously generated and then transferred to the GPU global RAM memory.

The outline of this paper is as follows. In Section \ref{sec:sabrmodel} the static and dynamic SABR models are described, with a new proposal for functional parameters in the dynamic case. In Section \ref{sec:calibration} the calibration to market data is detailed. In Section \ref{Monte_Carlo_GPU}, the Monte Carlo method and its parallel implementation in CUDA is described. In this work we use the SA method to calibrate the models. Depending on whether the cost function uses a direct expression or a Monte Carlo method, in Section \ref{sec:calibracionGpu} we discuss different techniques to parallelize the SA algorithm. In Section \ref{sec:numerical_results} we illustrate the performance of the implemented pricing technique for European call options. Next, we present results about calibration to real market data. Finally, the pricing of a cliquet option with the calibrated parameters is detailed.

%
%
%
%

\section{The SABR model}\label{sec:sabrmodel}

The SABR (\textit{Stochastic} $\alpha,\beta,\rho$) model was introduced in \cite{Hagan}, arguing that local volatility models could not reproduce market volatility smiles and that their predicted volatility dynamics contradicts market smiles and skews. The main advantage of the SABR model comes from its great simplicity compared to alternative stochastic volatility models \cite{Hagan}. The dynamics of the forward price and its volatility satisfy the system of stochastic differential equations
\begin{eqnarray}
dF_{t} \; =  \; \alpha_{t}F_{t}^{\beta}dW^{1}_t, && F_0= \hat{f}, \label{eq:model_staticSABR1}\\
d\alpha_{t} \; = \;  \nu \alpha_{t}dW^{2}_t, && \alpha_0=  \alpha, \label{eq:model_staticSABR2}
\end{eqnarray}
where $F_{t}=S_{t}e^{(r-y)(T-t)}$ denotes the \textit{forward} price of the underlying asset $S_t$, $r$ being the  constant interest rate and $y$ being the constant dividend yield. Moreover, $\alpha_{t}$ denotes the asset volatility process, $dW^{1}$ and $dW^{2}$ are two correlated Brownian motions with constant correlation coefficient $\rho$ (i.e. $dW_t^1 dW_t^2 = \rho dt$) and $S_0$ is the spot price of the asset. The parameters of the model are: $\alpha>0$ (the volatility's reference level), $0\leq\beta\leq1$ (the variance elasticity), $\nu>0$ (the volatility of the volatility) and $\rho$ (the correlation coefficient). Note the two special cases: $\beta=1$ (lognormal model) and $\beta=0$ (normal model).


\subsection{Static SABR model}
 The static SABR model corresponds to a constant parameters assumption. When working with options with the same maturity, the static SABR model provides good results \cite{Hagan}. The great advantage is that the following asymptotically approximated explicit formula for the implied Black-Scholes volatility can be obtained:
{\small
\begin{equation} \label{haganf}
\begin{split}
\sigma_{model}(K,\hat{f},T) =& \dfrac{\alpha}{(K\hat{f})^{(1-\beta)/2}\left[1 + \dfrac{(1-\beta)^{2}}{24}\ln^{2}\left(\dfrac{\hat{f}}{K}\right) +
 \dfrac{(1-\beta)^{4}}{1920}\ln^{4}\left(\dfrac{\hat{f}}{K}\right) + \cdots \right]} \cdot \left(\dfrac{z}{x(z)} \right) \cdot
\\
\\
& \left\{1 + \left[ \frac{(1-\beta)^{2}}{24}\frac{\alpha^{2}}{(K\hat{f})^{1-\beta}} + \frac{1}{4}
\frac{\rho \beta \nu \alpha}{(K\hat{f})^{(1-\beta)/2}} + \frac{2 - 3\rho^{2}}{24}\nu^{2}\right]\cdot T + \cdots\right\},
\end{split}
\end{equation}
}
where $z$  is a function of $K$, $\hat{f}$ and $T$ given by
\begin{equation} \label{zeta}
z = \frac{\nu}{\alpha}(K\hat{f})^{(1-\beta)/2}\ln \left(\dfrac{\hat{f}}{K}\right), \nonumber
\end{equation}
and
\begin{equation} \label{equis}
x(z) = \ln \left(\frac{\sqrt{1 - 2\rho z + z^{2}} + z - \rho}{1 - \rho}\right).
\end{equation}

In this work, we consider the following correction to (\ref{haganf}) proposed by Ob{\l}{\'o}j in \cite{obloj},
{\small
\begin{equation}
\label{eq:sigma-SABR}
\begin{split}
\sigma_{model}(K,\hat{f},T) =& \dfrac{1}{\left[1 + \dfrac{(1-\beta)^{2}}{24}\ln^{2}\left(\dfrac{\hat{f}}{K}\right) +
 \dfrac{(1-\beta)^{4}}{1920}\ln^{4}\left(\dfrac{\hat{f}}{K}\right) + \cdots \right]} \cdot \left(\dfrac{\nu\ln \left(\frac{\hat{f}}{K}\right)}{x(z)} \right) \cdot
\\
\\
& \left\{ 1 + \left[\frac{(1-\beta)^{2}}{24}\frac{\alpha^{2}}{(K\hat{f})^{1-\beta}} + \frac{1}{4}
\frac{\rho \beta \nu \alpha}{(K\hat{f})^{(1-\beta)/2}} + \frac{2 - 3\rho^{2}}{24}\nu^{2}\right]\cdot T + \cdots \right\},
\end{split}
\end{equation}
}
\noindent where the following new expression for $z$ is considered:
\begin{equation}
z = \frac{\nu \left(\hat{f}^{1-\beta}-K^{1-\beta}\right)}{\alpha(1-\beta)}, \nonumber
\end{equation}
and $x(z)$ is given by (\ref{equis}). The omitted terms after  $+ \cdots$ can be neglected, so that (\ref{eq:sigma-SABR}) turns to
\begin{equation}
\sigma_{model}(K,\hat{f},T) = \frac{1}{\omega}\left(1 + A_{1}\ln \left(\dfrac{K}{\hat{f}}\right) + A_{2}\ln^{2}\left(\dfrac{K}{\hat{f}}\right) + BT\right),
\label{vola_SABR}
\end{equation}
\noindent where the coefficients $A_{1}$, $A_{2}$ and $B$ are given by
\begin{eqnarray}
A_{1} &=& \displaystyle -\frac{1}{2}(1 - \beta - \rho \nu \omega), \nonumber\\
A_{2} &=& \displaystyle \frac{1}{12}\Big((1-\beta)^{2} + 3\big((1-\beta) - \rho \nu \omega \big) + \left(2 - 3\rho^{2}\right)\nu^{2}\omega^{2}\Big), \nonumber\\
B &=& \displaystyle \frac{(1-\beta)^{2}}{24}\frac{1}{\omega^{2}} + \frac{\beta \rho \nu}{4}\frac{1}{\omega} + \frac{2-3\rho^2}{24}\nu^{2}, \nonumber
\end{eqnarray}

\noindent and the value of $\omega$ is given by $\omega = \alpha^{-1} \hat{f}^{1-\beta}$.

\subsection{Dynamic SABR model and the choice of the functional parameters}\label{elec_func}

The main drawback of the static SABR model arises when market data for options with several maturities are considered. In this case, too large errors can appear. In order to overcome this problem, the following dynamic SABR model allows time dependency in some parameters \cite{Hagan}:
\begin{eqnarray}
       dF_{t} = \alpha_{t}F_{t}^{\beta}dW_t^1, && F_0=\hat{f} , \label{eq_SABRd1}\\
       d\alpha_{t} = \nu(t) \alpha_{t}dW_t^2, && \alpha_0 = \alpha, \label{eq_SABRd2}
\end{eqnarray}
where the correlation coefficient $\rho$ is also time dependent. As in the static SABR model, the dynamic one also provides the following expression to approximate the implied volatility \cite{Osajima-2007}:
{\small
\begin{equation}
\sigma_{model}(K,\hat{f},T) = \frac{1}{\omega}\left(1 + A_{1}(T)\ln\left(\frac{K}{\hat{f}}\right)
 + A_{2}(T)\ln^{2}\left(\frac{K}{\hat{f}}\right) + B(T)T\right),
\label{vola_SABRd}
\end{equation}
}
where
{\small
\begin{eqnarray}
A_{1}(T) &=& \displaystyle \frac{\beta - 1}{2} + \frac{\eta_{1}(T)\omega}{2}, \nonumber\\
A_{2}(T) &=& \displaystyle \frac{(1-\beta)^{2}}{12} + \frac{1-\beta -  \eta_{1}(T)\omega}{4} + \frac{4\nu_{1}^2(T) +3\left(\eta_{2}^2(T) - 3\eta_{1}^2(T)\right)}{24}\omega^{2}, \nonumber\\
B(T) &=& \displaystyle \frac{1}{\omega^2}\left(\frac{(1-\beta)^{2}}{24} + \frac{\omega \beta \eta_{1}(T)}{4} + \frac{2\nu_{2}^2(T) - 3\eta_{2}^2(T)}{24}\omega^{2}\right), \nonumber
\end{eqnarray}
}
with
{\small
\begin{equation}
\begin{array}{l}
\nu_{1}^2(T) = \displaystyle \frac{3}{T^{3}} \displaystyle \int_{0}^{T}{(T-t)^{2}\nu^{2}(t)dt},\quad
\nu_{2}^2(T) = \displaystyle \frac{6}{T^{3}} \displaystyle \int_{0}^{T}{(T-t)t\nu^{2}(t)dt},\\
\eta_{1}(T) = \displaystyle \frac{2}{T^{2}} \displaystyle \int_{0}^{T}{(T-t)\nu(t)\rho(t)dt},\quad
\eta_{2}^2(T) = \displaystyle \frac{12}{T^{4}} \displaystyle \int_{0}^{T}{\int_0^t{\left(\int_0^s{\nu(u)\rho(u)du}\right)^2ds}dt}.
\end{array}
\label{integrales_SABRd}
\end{equation}
}
 Note that if $\nu$ and $\rho$ are taken as constants, i.e. $\nu=\nu_0$ and $\rho=\rho_0$, then it follows that
$\nu_1(T) = \nu_2(T) = \nu_0$, $\eta_1(T) = \eta_2(T) = \nu_0\rho_0$ and the dynamic SABR model reduces to the static one.

The choice of the functions $\rho$ and $\nu$ in (\ref{integrales_SABRd}) constitutes a very important decision. The values of $\rho(t)$ and $\nu(t)$ have to be smaller for long terms ($t$ large) rather than for short terms ($t$ small). Thus, in this work we consider two possibilities with exponential decay:
\begin{itemize}
 \item {\bf Case I}: It is more classical and corresponds to the choice
\begin{equation}
\rho(t) = \rho_0e^{-at},\quad \nu(t) = \nu_0e^{-bt},
\label{el_fun_1}
\end{equation}
with $\rho_0 \in [-1,1]$, $\nu_0 > 0$, $a\geq 0$ and $b \geq 0$. In this case, the expressions of the functions $\nu_1^2$, $\nu_2^2$, $\eta_1$ and $\eta_2^2$, defined by (\ref{integrales_SABRd}), can be exactly calculated and are given by:
{\small
\begin{equation}
\begin{array}{l}
\nu_1^2(T) = \displaystyle \frac{6\nu_0^2}{(2bT)^3}\left[\left((2bT)^2/2 - 2bT + 1 \right) - e^{-2bT}\right],\\
\\
\nu_2^2(T) = \displaystyle \frac{6\nu_0^2}{(2bT)^3}\left[2(e^{-2bT} - 1) + 2bT(e^{-2bT} + 1)\right],\\
\\
\eta_1(T) = \displaystyle \frac{2\nu_0\rho_0}{T^2(a+b)^2}\left[e^{-(a+b)T} - \big(1 - (a+b)T\big)\right],\\
\\
\eta_2^2(T) = \displaystyle \frac{3\nu_0^2\rho_0^2}{T^4(a+b)^4}\left[e^{-2(a+b)T}-8 e^{-(a+b)T}+\Big(7+2(a+b)T\big(-3+(a+b)T\big)\Big)\right].
\end{array}
\label{eq:exact-integral-expressions}
\end{equation}
}
 \item {\bf Case II}: In the present paper we propose the original and more general choice
\begin{equation}
\rho(t) = (\rho_0 + q_{\rho}t)e^{-at} + d_{\rho},\quad
\nu(t) = (\nu_0 + q_{\nu}t)e^{-bt} + d_{\nu}.
\label{el_fun_3}
\end{equation}

In this case, the symbolic software package Mathematica allows to calculate exactly the functions $\nu_1^2$, $\nu_2^2$ and $\eta_1$ (see Appendix \ref{anexo_funciones_casogeneral}). However, an explicit expression for  $\eta_2^2$ cannot be obtained, and therefore we use an appropriate quadrature formula for its approximation.

In order to guarantee that the correlation $\rho(t) \in [-1,1]$ and the volatility $\nu(t)>0$ for the involved parameters, an adequate optimization algorithm has to be used during calibration.
\end{itemize}

\section{Calibration of the SABR model} \label{sec:calibration}

The goal of calibration is to fit the model parameters to reproduce the market prices or volatilities as close as possible. In order to calibrate the model, we choose a set of vanilla options on the same underlying asset and the market prices are collected at the same moment. We can either consider only one maturity or several maturities. In the second case, a dynamic SABR model should be applied. Next, the model is used to price other options (such as exotic options). The computational cost increases with the increasing complexity of the pricing models.

Hereafter we denote by $Data_{market}(K_j,\hat{f},T_i)$ the observed market data, for the maturity $T_i$ ($i=1,\ldots, n$) and the strike $K_j$ ($j=1,\ldots, m_i$), where $n$ denotes the number of maturities and $m_i$ the number of strikes for the maturity $T_i$. The market implied volatility is denoted by $\sigma_{market}(K_j,\hat{f},T_i)$ and the corresponding market option price by $V_{market}(K_j,\hat{f},T_i)$. Moreover, we denote by $Data_{model}(K_j,\hat{f},T_i)$ the model value ($V_{model}(K_j,\hat{f},T_i)$ or $\sigma_{model}(K_j,\hat{f},T_i)$) for the same option.

The calibration process tries to obtain a set of model parameters that minimizes the error between market and model values for a given error measure. In order to achieve this target we must follow several steps:
\begin{itemize}
	\item Decide how to perform the calibration process, in prices or in volatilities.
	\item Choose market data that should be highly representative of the market situation.
	\item Decide which error measure will be used to compare model and market values:
   \begin{itemize}

\item If the calibration is made for one maturity $T_i$, we use the cost function
%
%
%
{\small
    	\begin{equation}
	f_{i,E}(\pmb x)=\displaystyle  \sum_{j=1}^{m_i} \left( \frac{ Data_{market}(K_j,\hat{f},T_i) - Data_{model}(K_j,\hat{f},T_i)}{Data_{market}(K_j,\hat{f},T_i)}  \right)^2(\pmb x),
	\label{eq:error_calibrate_individual_relativo}
	\end{equation}
}
   \noindent where $\pmb x$ denotes the parameters to calibrate.
\item If the calibration is made for a set of maturity dates the cost function we use is
{\small
	\begin{equation}
	f_E(\pmb x)=\displaystyle \sum_{i=1}^{n} \sum_{j=1}^{m_i} \left( \frac{ Data_{market}(K_j,\hat{f},T_i) - Data_{model}(K_j,\hat{f},T_i)}{  Data_{market}(K_j,\hat{f},T_i) }\right)^2(\pmb x).
	\label{eq:error_calibrate_relativo}
	\end{equation}
}
   \end{itemize}

	\item Choose the (local or global) optimization algorithm to minimize the error.
	\item Fix (if it is convenient) some of the parameters on beforehand, by taking into account the previous ex\-pe\-rien\-ce or the existing information.
	\item Calibrate and compare the obtained results. If they are satisfactory, the parameters are accepted and used for pricing more complex financial instruments.
\end{itemize}

An advantage of the SABR model is the existence of an asymptotic approximation formula for the implied volatility that can be used in the calibration. It is important to take into account the meaning of the different model parameters \cite{Hagan,west}. The value of $\beta$ is related to the type of the underlying stochastic process of the model and it is usually fixed on beforehand. For example, $\beta=1$ (lognormal model) is mostly used in equity and currency markets, like foreign exchange markets (for example, EUR/USD). The choice $\beta=0.5$ corresponds to a CIR model and is used in US interest rate desks. The value $\beta=0$ (normal model) is commonly used to manage risks, as for example in Yen interests rate markets. Alternatively, $\beta$  can be computed from historical data of the at-the-money volatilities, $\sigma_{ATM}=\sigma_{model}(\hat{f},\hat{f},T)$, that can be obtained by taking $K=\hat{f}$ in \eqref{eq:sigma-SABR}. Next, after some computations the identity
$$ \ln \sigma_{ATM}=\ln \alpha -(1-\beta)\ln \hat{f},$$
is obtained and $\beta$ can be computed from a log-log regression of $\hat{f}$ and $\sigma_{ATM}$ (see \cite{west}, for details). Once $\beta$ has been fixed and the term $B T$ in \eqref{eq:sigma-SABR} has been neglected, the value of $\alpha$ can be computed by using the value of $\sigma$ at-the-money ($\sigma_{ATM} = \sigma_{model}(\hat{f},\hat{f},T)$) with the following approximation: $$ \alpha\approx \hat{f}^{(1-\beta)} \sigma_{model}(\hat{f},\hat{f},T). $$
In the case $\beta=1$,  $\alpha$ is  equal to the at-the-money volatility \cite{Hagan}. Thus, by fixing $\beta=1$ and  $\alpha$, the (possibly large) calibration cost is reduced: only $\rho$ and $\nu$ are calibrated.

 In this work the whole set of SABR parameters is calibrated, thus leading to a more time consuming calibration. Nevertheless, the GPUs technology highly speeds up this procedure.

\section{Pricing with Monte Carlo using GPUs}\label{Monte_Carlo_GPU}

Monte Carlo technique for the SABR model involves the simulation of a huge number of forward and volatility paths. For a given option data and a given initial value of the forward $F_0=\hat{f}=S_0 e^{(r-y)(T-t)}$, the European option price at the time $t=0$ is equal to $V(S_0,K) = e^{-r T}\mathbb{E}\big(V(S_T,K)\big)$, where $V(S_T,K)$ is the option payoff (for example, $V(S_T,K) = \max\{S_T-K,0\}$ for a European call option).

In order to discretize the stochastic differential equations of the SABR model, we take a constant time step $\Delta t$, such that $M=\frac{T}{\Delta t}$ denotes the number of subintervals in $[0,T]$. In order to preserve some properties of the continuous model (such as positivity), we use the following log-Euler discretization:
\begin{eqnarray}
       \alpha_{i+1} & = & \alpha_ie^{\nu(t_i) Z^1_i \sqrt{\Delta t} - \nu^2(t_i)\Delta t/2}, \nonumber\\
       \hat{\nu} & = & \alpha_iF_i^{\beta-1}, \label{dis_SABRd} \\
       F_{i+1} & = & F_ie^{\hat{\nu} \left(\rho(t_i)Z^{1}_i+ Z^{2}_i\sqrt{1-\rho^2(t_i)}\right)\sqrt{\Delta t} - \hat{\nu}^2 \Delta t/2}, \nonumber
\end{eqnarray}
where $Z^{1}_i$ and $Z^2_i$ denote two independent standard normal distribution samples. In \cite{ChenOsterleeHans_2012_lowbias} it is pointed out that the log-Euler scheme may become unstable for specific time-steps. So, the simulation of the conditional distribution of the SABR model over a discrete time step (which is proved to be a transformed squared Bessel process) and the use of an approximation formula for the integrated variance are proposed. However, we did not find this unstable behavior in our experiments. In case of unstability, this alternative low-bias simulation scheme in \cite{ChenOsterleeHans_2012_lowbias} could be adapted.

\noindent Note that $S_T=F_T$. If $N$ denotes the number of simulated paths, then the option price is given by
\begin{equation} \label{payoff}
\widehat{V}(S_0,K) = \displaystyle e^{-rT}\frac{1}{N}\sum_{j=1}^N V(S_{T}^j,K).
\end{equation}

\noindent Option pricing with the Monte Carlo method is an accurate and robust technique. However, as Monte Carlo method exhibits an order of convergence equal to $1/\sqrt{N}$, a large number of paths is required to obtain precise results, therefore leading to not affordable computational costs in real practice. For this reason, option pricing with Monte Carlo is usually implemented in high performance systems. Here, we propose a parallel version of Monte Carlo efficiently implemented in CUDA by using the XORWOW pseudo-random number generator (xorshift RNG) included in the Nvidia CURAND library \cite{art:MarsagliaXorshift_03}.



The process can be summarized as follows:

\begin{enumerate}
 \item Firstly we read the input data and send it to constant memory: Monte Carlo parameters, SABR parameters and market data. Next, in order to store the computed final price of each simulation, we allocate a global memory vector of size $N$.


\item We compute a Monte Carlo kernel for the generation of the different paths, where uniform random numbers are generated ``on the fly'' inside this kernel, calling the {\tt curand\_uniform()} CURAND function. This allows random numbers to be generated and used by the Monte Carlo kernel without requiring them to be written to and then read from global memory. Then Box-Muller method is used to transform the random uniform distributed numbers in normally distributed ones.

\item The option price is computed with the {\tt thrust::transform\_reduce} method of the Thrust
Library \cite{ref:thrust} (included in the Nvidia toolkit since its version 4.0). With
this function we apply kernel fusion reduction kernels. Then, the sum in (\ref{payoff}) is computed by a standard \texttt{plus} reduction also available in the Thrust Library. By fusing the payoff operation with the reduction kernel we have a highly optimized implementation which offers the same performance as hand-written kernels. All these operations are performed inside the GPU, so that transfers to the CPU memory are avoided.


\end{enumerate}

\section{Calibration of the parameters using GPUs}\label{sec:calibracionGpu}

The calibration of the SABR model parameters can be done using the implied volatility formula or the Monte Carlo simulation method. Usually, in trading environments the second one is not used, mainly due to its high execution times. However, if we have a parallel and efficient implementation of the Monte Carlo method, we can consider its usage in the calibration procedure.

In this work, the calibration of the parameters has been done with a Simulated Annealing (SA) stochastic global optimization method. The SA algorithm consists of a temperature loop (see \cite{BM-95} for details), where the equilibrium state at each temperature is computed using the Metropolis algorithm.

In this paper, we propose the two calibration techniques described in next paragraphs.

\subsection{Calibration with Technique I}
 In this calibration technique we use the implied volatility formula and an efficient implementation in GPU of the Simulated Annealing algorithm. Thus, we consider $Data_{market}=\sigma_{market}$ and $Data_{model}=\sigma_{model}$ in the cost functions \eqref{eq:error_calibrate_individual_relativo} and \eqref{eq:error_calibrate_relativo}, respectively, both for individual and joint calibrations. Note that $\sigma_{market}$ denotes the market quoted implied volatilities and $\sigma_{model}$ the model ones, \eqref{vola_SABR} and \eqref{vola_SABRd} for the static and dynamic models, respectively.

As it is well known in the literature, SA involves a great computational cost. For this reason, the calibration of the pa\-ra\-me\-ters has been done using the CUSIMANN library \cite{page:cusimann}, that is a parallel implementation of the SA algorithm proposed in \cite{art:SA-GPU}; specifically the so called synchronous implementation. In \cite{art:SA-GPU}, authors have analyzed different parallel SA implementations and have concluded that the synchronous version is the most appropriated for the GPU technology.

The idea of this synchronous parallel version of SA is the following: the algorithm starts from a fixed initial solution $\pmb x_0$, so that each thread runs independently a Markov chain of constant length $L$ until reaching the next level of temperature. As the temperature is fixed, each thread actually performs a Metropolis process. Once all threads have finished, they report their corresponding final states $\pmb x^p$ and the value $f(\pmb x^p)$, $p=0,\ldots, w-1$ (where $w$ denotes the number of threads). Next, a reduce operation to obtain the minimum of the cost function is performed. So, if the minimum is obtained at a particular thread $p^{\star}$ then $\pmb x^{p^{\star}}$ is used as starting point for all threads at the following temperature level. In \cite{art:SA-GPU} the implementation on CUDA of this parallel version is detailed.


A host can have several discrete graphic cards attached. When dealing with computationally expensive processes, as in this case, it is interesting to use the computing power of these extra GPUs to further reduce the execution times. One possible approach is to use OpenMP \cite{ref:openmp} to control the GPUs inside one node. The strategy is to launch as many CPU threads as GPUs available on a node. Thus, each CPU thread handles a GPU.

\begin{figure}[!htb]
\begin{center}
\includegraphics[height=5cm]{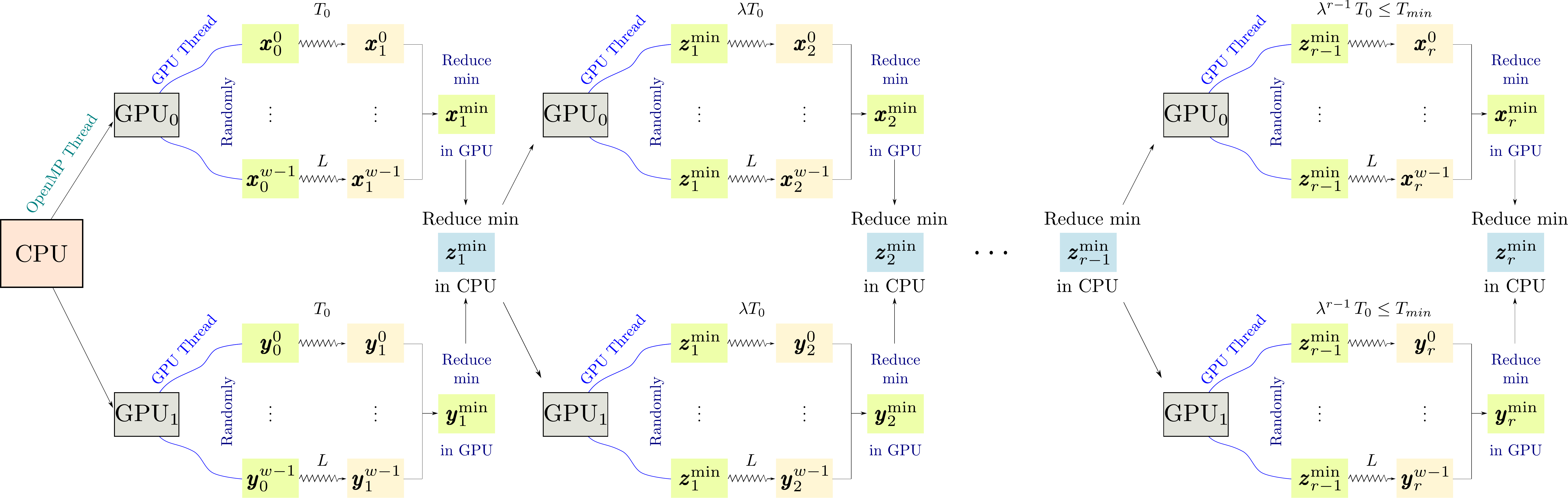}
\caption{Sketch of the parallel SA algorithm using two GPUs and OpenMP.}
\label{fig:sketchParaSincrono_2GPUs}
\end{center}
\end{figure}

The idea of this new synchronous parallel version of SA considering OpenMP and several GPUs is the following; by simplicity we will consider two GPUs (see Figure \ref{fig:sketchParaSincrono_2GPUs}): in each GPU the algorithm starts from a fixed initial point, and each thread runs independently a Markov chain of constant length $L$ until reaching the next level of temperature. As the temperature is fixed, each thread actually performs a Metropolis process. Once all the threads have finished, in each GPU a reduction operation is made. Then, each GPU reports the corresponding final state $\pmb x_{*}^{\text{min}}$ and $\pmb y_{*}^{\text{min}}$, and in the CPU another reduction operation is done to get the minimum point $\pmb z_{*}^{\text{min}}$. This point is used as starting point for all GPUs threads at the following temperature level.

 Note that if we had a cluster of GPUs, we could exploit an additional level of parallelism. The final configuration which merges all the ideas is to use the MPI \cite{ref:mpi} for intercommunication between the nodes of a GPU cluster and while to use OpenMP inside the node to control several graphic cards. When we double the number of GPUs used, the benefit of this approach is that the execution times are reduced by about a half, or, from another point of view, we could double the number of function evaluations ``without'' increasing the runtime. Indeed, the time reduction would be somewhat less than a half because this extra parallelization involves a small additional cost.

\subsection{Calibration with Technique II}
In this calibration technique the cost function is computed in GPU by a Monte Carlo method. It turns out to be necessary when an expression for the implied volatility is not available, as in the dynamic SABR model when considering the {\bf Case II} (see Section \ref{elec_func}). For the joint calibration, in the cost function \eqref{eq:error_calibrate_relativo} we consider $Data_{market}=V_{market}$ and $Data_{model}=V_{model}$; $V_{market}$ being the option market price and $V_{model}$ being the computed price with the Monte Carlo parallel implementation (Section \ref{Monte_Carlo_GPU}).

As the Monte Carlo method is carried out inside the GPU and Fermi GPUs do not allow nesting kernels, the SA minimization algorithm has to be run on CPU. In order to use all available GPUs in the system, we propose a CPU SA parallelization using OpenMP. So, each OpenMP SA thread uses a GPU to assess on the Monte Carlo objective function. This approach can be easily extrapolated to a cluster of GPUs using, for example, MPI. Particularly, in this case the SA procedure rejects all points violating the constraints over $\rho(t)$ and $\nu(t)$, i.e. $\rho(t) \in [-1,1]$ and $\nu(t) > 0$, $\forall t$ $0<t\leq T$.


\section{Numerical results}\label{sec:numerical_results}
From now on we denote by {\tt SSabr} and {\tt DSabr} the static and dynamic SABR models, respectively. As detailed in Section \ref{elec_func}  concerning the dynamic SABR model, depending on the choice of the time dependent functions $\rho$ and $\nu$, we denote by:
\begin{itemize}
 \item {\tt DSabr\_I}: dynamic SABR model in the {\bf Case I}, where $\rho$ and $\nu$ are given by \eqref{el_fun_1}.
\item {\tt DSabr\_II}: dynamic SABR model in the {\bf Case II} (general case), where $\rho$ and $\nu$ are given by \eqref{el_fun_3}.
\end{itemize}

Moreover, depending on the technique used to calibrate the model (see Section \ref{sec:calibracionGpu}), we use the notation:
\begin{itemize}
 \item {\tt T\_I}: {\bf Technique I}, where the model parameters are calibrated with the implied volatility formula and using the efficient implementation of Simulated Annealing in GPU.
\item {\tt T\_II}: {\bf Technique II}, where the cost function is evaluated by the efficient implementation in GPUs of the Monte Carlo pricing method, as detailed in Section \ref{Monte_Carlo_GPU}.
\end{itemize}

Numerical experiments have been performed with the following hardware and software configurations: GPUs Nvidia Geforce GTX470, a CPU Xeon E5620 clocked at 2.4 Ghz with 16 GB of RAM, CentOS Linux, Nvidia CUDA SDK 4.0 and GNU C/C++ compilers 4.1.2.

\subsection{Pricing European options}

In this section we focus on pricing European options with both SABR models, when using the GPU Monte Carlo method.
The fixed data are $S_0=2257.37$, $K=2257.37$, $r=0.018196$, $y=0.034516$, $T=0.495890$, while the SABR parameters are:
\begin{itemize}
 \item {\tt SSabr}: $\alpha= 0.375162$, $\beta=0.999999$, $\nu=0.331441$ and $\rho=-0.999999$.
 \item {\tt DSabr\_I}: $\alpha=0.393329$, $\beta=1.0$, $\nu_0=0.941565$, $\rho_0=-1.0$, $a=0.001$ and $b=1.246906$.
 \item {\tt DSabr\_II}: $\alpha=0.398436$, $\beta=0.999579$, $\nu_0=1.285129$, $\rho_0=-0.964678$, $a=0.0$, $b=2.059560$, $d_{\rho}=0.101632$, $d_{\nu}=-0.086294$, $q_{\rho}=0.0$ and $q_{\nu}=1.302296$.
\end{itemize}

In Table \ref{tab:one_timesteps} the results for the three presented models are shown. They correspond to $2^{20}$ simulations, $\Delta t=1/250$ ($123$ time steps) and both single and double precision computations. In Table \ref{tab:speepup_one} the execution times for CPU an GPU programs are compared. In single precision, the maximum speedup varies, approximately, from $357$ to $567$ times depending on the model. In double precision, the maximum speedup is around of $74$ times in all models. Table \ref{tab:oneoption_increasing_strikes} shows the variation of execution times when increasing the number of strikes, considering the {\tt DSabr\_II} model with $\Delta t= 1/250$. Note that execution times hardly vary: pricing $41$ options by generating $2^{20}$ paths takes around $0.28$ and $0.86$ seconds in single and double precision, respectively, while pricing one option takes $0.25$ and $0.84$ seconds. Therefore, Monte Carlo pricing method results to be affordable for calibration.

\begin{table}[!htb]
\centering
{\footnotesize
\begin{tabular}{|c || c|c || c|c || c|c|}
\hline
  Precision& \multicolumn{2}{|c||}{\tt SSabr} & \multicolumn{2}{|c||}{\tt DSabr\_I} & \multicolumn{2}{|c|}{\tt DSabr\_II} \\
\hline\hline
 \, &Result & {\tt RE}& Result & {\tt RE} & Result & {\tt RE}  \\
\hline
Single &$224.544601$ & $0.005944$ & $222.432648$ & $0.015294$ & $224.652390$ & $0.005467$  \\
\hline
 Double &$224.545954$ & $0.005938$ & $222.434009$ & $0.015288$ & $224.653642$ & $0.005461$  \\
\hline
\end{tabular}
}
\caption{Pricing results for European options. {\tt RE} denotes the relative error with respect to the reference value $225.887329$.}
\label{tab:one_timesteps}
\end{table}

\begin{table}[!htb]
\centering
{\footnotesize
\begin{tabular}{|c || c|c|c || c|c|c|}
\hline
\multicolumn{7}{|c|}{\tt SSabr}  \\
\hline
Paths & \multicolumn{3}{|c||}{Single} & \multicolumn{3}{|c|}{Double} \\
\hline
\, & CPU & GPU & Speedup & CPU & GPU & Speedup \\
\hline
$2^{16}$ & $2.890$ & $0.142$ & $\times 20.313$ & $2.762$ & $0.175$ & $\times 15.783$ \\
\hline
$2^{20}$ & $46.287$ & $0.206$ & $\times 223.726$ & $44.142$ & $0.723$ & $\times 61.054$ \\
\hline
$2^{24}$ & $721.592$ & $1.271$ & $\times 567.373$ & $704.617$ & $9.474$ & $\times 74.373$ \\
\hline
\hline
\multicolumn{7}{|c|}{\tt DSabr\_I}\\
\hline
Paths & \multicolumn{3}{|c||}{Single} & \multicolumn{3}{|c|}{Double} \\
\hline
\, & CPU & GPU & Speedup & CPU & GPU & Speedup \\
 \hline
$2^{16}$ &  $3.638$ & $0.144$ & $\times 25.112$ & $3.358$ & $0.179$ & $\times 18.759$ \\
\hline
$2^{20}$ &  $58.228$ & $0.255$ & $\times 227.674$ & $53.663$ & $0.840$ & $\times 63.884$ \\
\hline
$2^{24}$ &  $929.254$ & $2.047$ & $\times 453.918$ & $860.556$ & $11.515$ & $\times 74.733$ \\
\hline\hline
\multicolumn{7}{|c|}{\tt DSabr\_II}\\
\hline
Paths & \multicolumn{3}{|c||}{Single} & \multicolumn{3}{|c|}{Double} \\
\hline
\, & CPU & GPU & Speedup & CPU & GPU & Speedup \\
\hline
$2^{16}$ & $3.694$ & $0.145$ & $\times 25.388$ & $3.392$ & $0.182$ & $\times 18.637$ \\
\hline
$2^{20}$ & $59.790$ & $0.287$ & $\times 207.762$ & $54.144$ & $0.853$ & $\times 64.474$ \\
\hline
$2^{24}$ & $949.459$ & $2.658$ & $\times 357.095$ & $868.156$ & $11.590$ & $\times 74.905$ \\
\hline
\end{tabular}
}
\caption{Pricing European options. Execution times for CPU and GPU calibration (in seconds), con\-si\-de\-ring single and double precision with $\Delta t= 1/250$.}
\label{tab:speepup_one}
\end{table}


\begin{table}[!htb]
\centering
{\footnotesize
\begin{tabular}{|c ||c| c||  c| c||  c| c||}
\hline
\multirow{2}{*}{Num. Strikes} & \multicolumn{2}{|c||}{$2^{16}$} & \multicolumn{2}{|c||}{$2^{20}$} & \multicolumn{2}{|c||}{$2^{24}$}\\
 & Single & Double & Single & Double& Single & Double\\
\hline
 $1$  & $0.134374$ & $0.183608$ & $0.259048$ & $0.841376$& $2.234662$ & $11.595257$\\
 \hline
$5$ & $0.135871$ & $0.184051$ & $0.262661$ & $0.850931$& $2.242463$ & $11.605406$\\
 \hline
 $41$  & $0.157753$ & $0.190635$ &$0.280509$ & $0.862674$  & $2.265047$ & $11.645788$\\
\hline
\end{tabular}
}
\caption{Pricing with {\tt DSabr\_II} model in GPU. Influence of the number of strikes in computational times (time in seconds), $\Delta t=1/250$. We consider $1$, $5$ and $41$ strikes, with values $K$ (in \% of $S_0$) of $\{100\}$, $\{96,\,98, 100,\,102,\,104\}$ and $\{80,\,81\, \ldots,\,119,\,120\}$, respectively.}
\label{tab:oneoption_increasing_strikes}
\end{table}

\subsection{Calibration}

In order to check the accuracy and performance of the GPU calibration code, the calibration of the SABR  model to the volatility surfaces generated by the EURO STOXX 50 index and EUR/USD foreign exchange rate has been tested. For both data, we detail the calibrated parameters for {\tt SSabr}, {\tt DSabr\_I} and {\tt DSabr\_II} models, the computational times and the comparisons between the market volatilities or prices and those ones with the model after parameter calibration.

\subsubsection{EURO STOXX 50 index}\label{subsub_eurostoxx}


The data corresponds to EURO STOXX 50 quotes of December of 2011. The asset spot value is $2311.1$ \euro. In Tables \ref{eurostoxx1} and \ref{eurostoxx1_vola} of the Appendix \ref{Append:MarketData}, the interest rates, dividend yields and implied volatilities for $3$, $6$, $12$ and $24$ months maturities are shown. Next, we present the calibrated parameters for these data.


\begin{enumerate}
\item {\bf Calibration of the {\tt SSabr} model using the technique {\tt T\_I}}

By using the technique {\tt T\_I} and the asymptotic expression (\ref{vola_SABR}) in the cost function, the individual calibration (all strikes and one maturity) of the {\tt SSabr} model has been carried out. In Table \ref{tab:staticCallibrateDataEURO STOXX 50_Individual} the calibrated parameters are detailed. For each maturity ($3$, $6$, $12$ and $24$ months), Figure \ref{fig:staticIndividualStoxx} shows market and model volatilities with the parameters of Table \ref{tab:staticCallibrateDataEURO STOXX 50_Individual}.

\begin{table}[!htb]
\centering
{\footnotesize
\begin{tabular}{|c|c|c|c|c|}
\hline
\, & $3$ months & $6$ months & $12$ months & $24$ months\\
\hline\hline
 $\alpha$ & $0.298999$ & $0.302060$ & $0.289271$ & $0.277844$\\
\hline
$\beta$ & $1.0$& $1.0$ & $1.0$ & $1.0$\\
\hline
$\nu$ & $0.382558$ & $0.381724$ & $0.308560$ & $0.264178$\\
\hline
$\rho$&$-1.0$& $-1.0$ & $-0.999729$ & $-1.0$ \\
\hline
\end{tabular}
}
\caption{EURO STOXX 50. {\tt SSabr} model: Calibrated parameters for each maturity.}
\label{tab:staticCallibrateDataEURO STOXX 50_Individual}
\end{table}

\begin{figure}[h!]
\centering
\subfigure{\includegraphics[height=5cm]{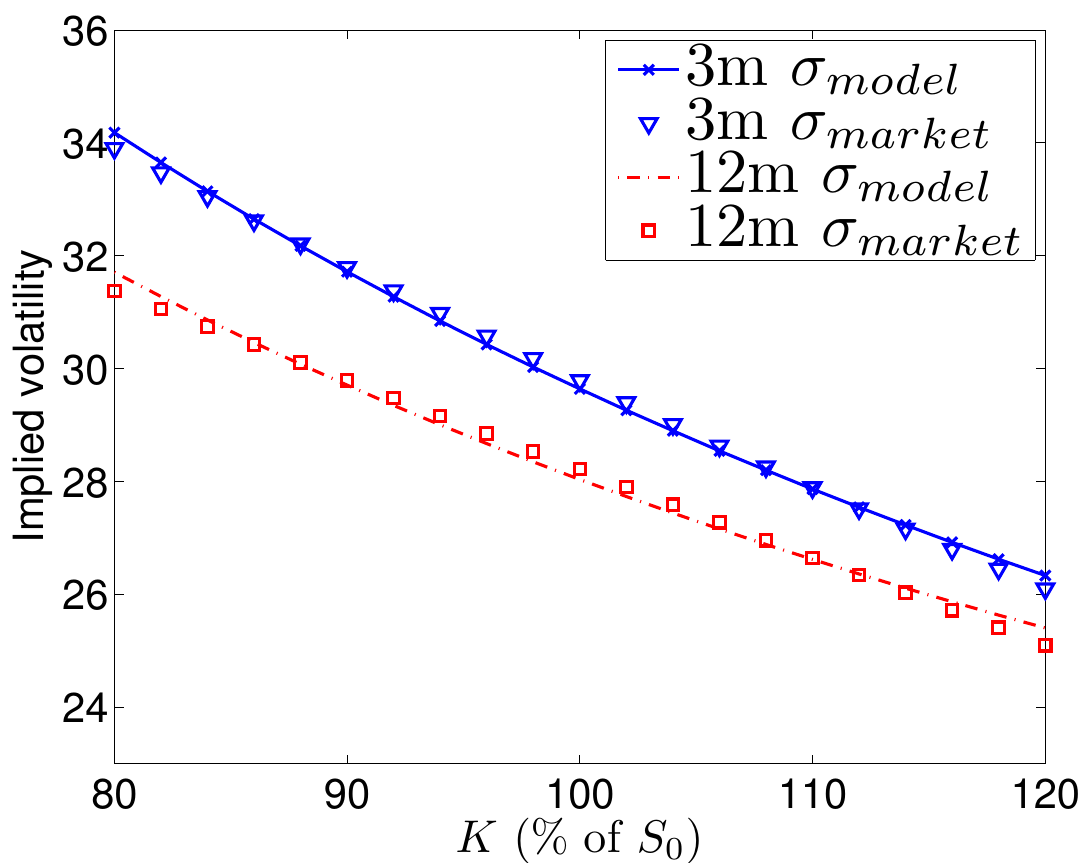}}
\subfigure{\includegraphics[height=5cm]{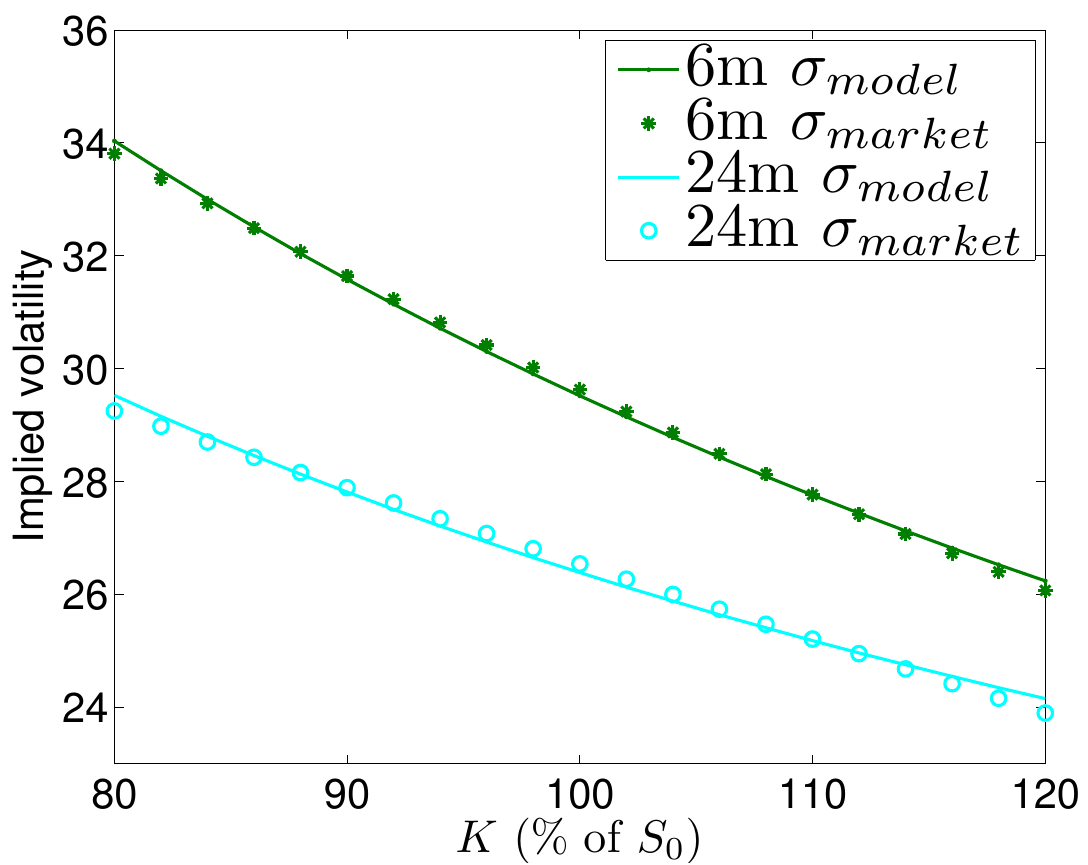}}
\caption{EURO STOXX 50. {\tt SSabr} model: $\sigma_{model}$ vs. $\sigma_{market}$ for the whole volatility surface. Maturities: 3 and 12 months (left), 6 and 24 months (right).}
\label{fig:staticIndividualStoxx}
\end{figure}

Table \ref{tab:static_24_Stoxx} shows the computational times and the speedups of the calibration process for $T=24$ months. As expected, the speedup is near $7.7$ when using OpenMP with $8$ threads, while the speedup respect to the CPU code is around $190$ when we use one GPU. Furthermore, if we use two GPUs then the computational time is additionally almost divided by two.

\begin{table}[!htb]
\centering
{\footnotesize
\begin{tabular}{|c|c|c|c|c|c|c|}
\hline
& CPU ($1$ thread) &$2$ threads &$4$ threads & $8$ threads &GPU  & $2$ GPUs\\
\hline
 Time ($s$) &$4172.746$ & $2096.857$& $1056.391$& $545.456$& $21.955$ & $12.609$\\
\hline
Speedup  & -  & $1.99$ & $3.95$ & $7.65$& $190.06$& $330.93$\\
\hline
\end{tabular}
}
\caption{EURO STOXX 50. {\tt SSabr} model: Performance of OpenMP vs. GPU versions, in single precision for $T = 24$ months.}
\label{tab:static_24_Stoxx}
\end{table}

\item {\bf Calibration of the {\tt DSabr\_I} model using the technique {\tt T\_I}}

The joint calibration for all strikes and maturities, using the {\tt DSabr\_I} model is made with the GPU version of SA and using the asymptotic formula \eqref{vola_SABRd},  when $\rho$ and $\nu$ are given by (\ref{el_fun_1}). In Table \ref{tab:dynamicSetCallibrateDataEURO STOXX 50} the calibrated parameters are detailed. In Figure \ref{fig:dynamicAllStoxx}, the whole volatility surface for all maturities is shown. By using the parameters in Table \ref{tab:dynamicSetCallibrateDataEURO STOXX 50}, for several strikes. Table \ref{tab:dynamicAllStoxx} shows the market volatilities ($\sigma_{market}$) vs. the model ones ($\sigma_{model}$, computed with formulas \eqref{vola_SABRd}, \eqref{el_fun_1} and \eqref{eq:exact-integral-expressions}). For single precision, the maximum relative error is $7.608205\times 10^{-2}$ and the mean relative error is $2.073025\times 10^{-2}$.

\begin{table}[!htb]
\centering
{\footnotesize
\begin{tabular}{|c|c |c|c|c|c|}
\hline
$\alpha=0.294722$ & $\beta=1.0$ & $\rho_0=-1.0$ & $\nu_0=0.388539$ &$a=0.001000$ & $b=0.131466$\\
\hline
\end{tabular}
}
\caption{EURO STOXX 50. {\tt DSabr\_I} model: Calibrated parameters.}
\label{tab:dynamicSetCallibrateDataEURO STOXX 50}
\end{table}

\begin{figure}[h!]
\centering
\subfigure{\includegraphics[height=5cm]{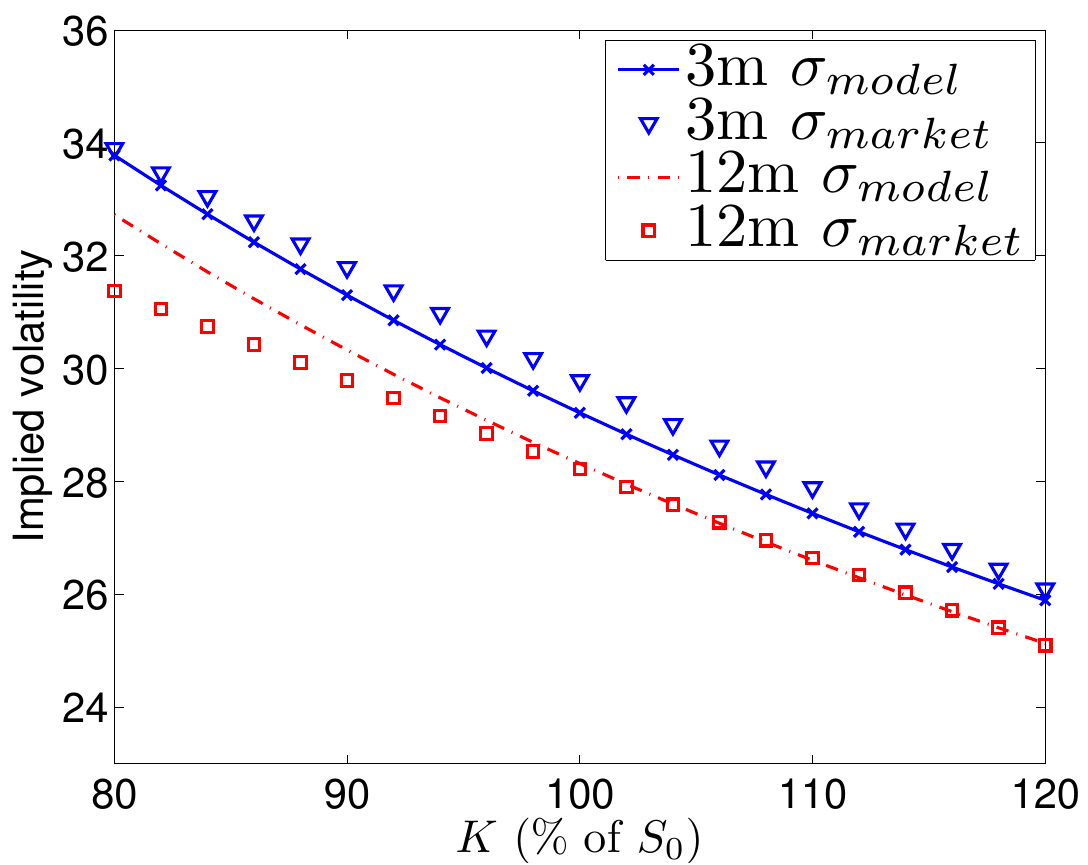}}
\subfigure{\includegraphics[height=5cm]{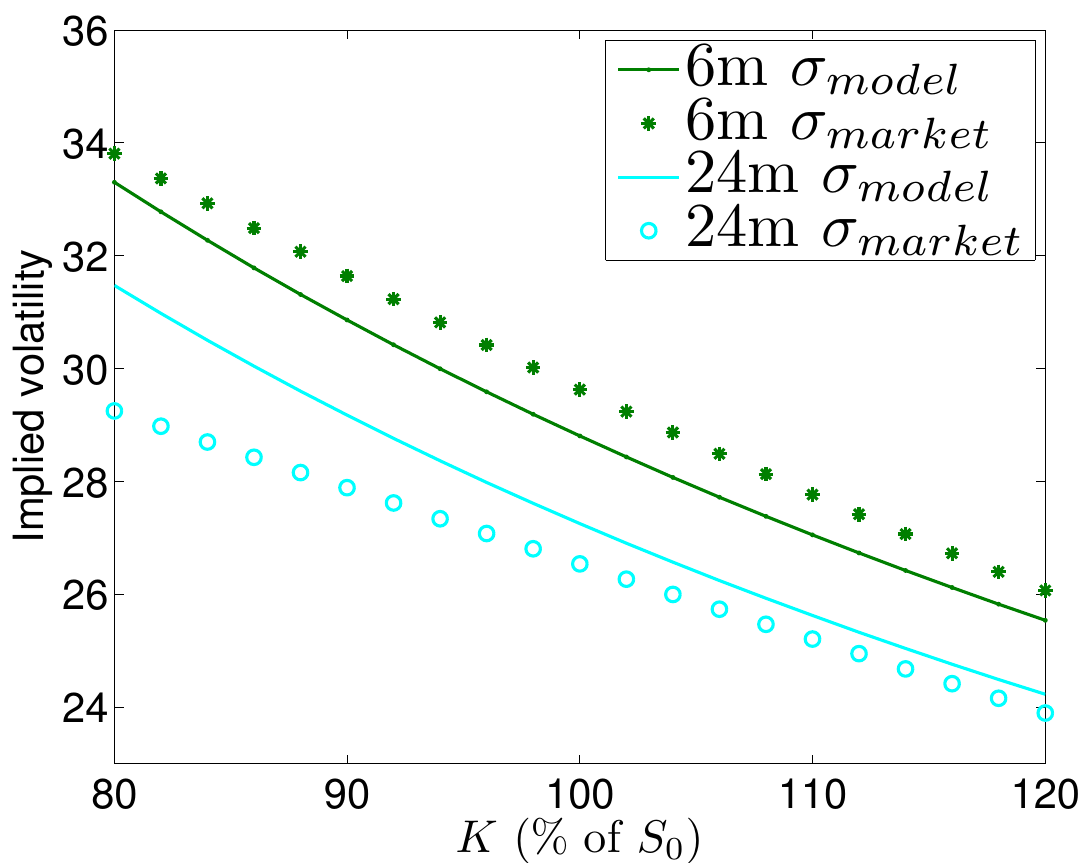}}
\caption{EURO STOXX 50. {\tt DSabr\_I} model: $\sigma_{model}$ vs. $\sigma_{market}$ for the whole volatility surface. Maturities: 3 and 12 months (left), 6 and 24 months (right).}
\label{fig:dynamicAllStoxx}
\end{figure}

\begin{table}[!htb]
\centering
{\footnotesize
\begin{tabular}{|c|| c|c|c ||c|c|c |}
 \hline
$K$ (\% of $S_0$) & \multicolumn{3}{|c|}{$3$ months }& \multicolumn{3}{|c|}{$6$ months } \\
\hline
\, & $\sigma_{market}$ & $\sigma_{model}$ & $\frac{|\sigma_{market}-\sigma_{model}|}{\sigma_{market}}$
         & $\sigma_{market}$ & $\sigma_{model}$ &$\frac{|\sigma_{market}-\sigma_{model}|}{\sigma_{market}}$ \\
\hline
$88\%$ & $32.21$ & $31.7628$ & $1.388389\times 10^{-02}$ & $32.07$ & $31.3150$ & $2.354225\times 10^{-02}$ \\
\hline
$100\%$ & $29.79$ & $29.2166$ & $1.924807\times 10^{-02}$ & $29.63$ & $28.8068$ & $2.778265 \times 10^{-02}$ \\
\hline
$112\%$ & $27.52$ & $27.1094$ & $1.492006\times 10^{-02}$ & $27.42$ & $26.7345$ & $2.500000 \times 10^{-02}$ \\
\hline
\hline
$K$ (\% of $S_0$) & \multicolumn{3}{|c|}{$12$ months } & \multicolumn{3}{|c|}{$24$ months}\\
\hline
\,& $\sigma_{market}$ & $\sigma_{model}$ & $\frac{|\sigma_{market}-\sigma_{model}|}{\sigma_{market}}$
       & $\sigma_{market}$ & $\sigma_{model}$ & $\frac{|\sigma_{market}-\sigma_{model}|}{\sigma_{market}}$\\
\hline
$88\%$ & $30.11$ & $30.7756$ & $2.210561\times 10^{-02}$ & $28.16$ & $29.6026$ & $5.122869\times 10^{-02}$ \\
\hline
$100\%$ & $28.22$ & $28.3187$ & $3.497519\times 10^{-03}$ & $26.54$ & $27.2549$ & $2.693670\times 10^{-02}$ \\
\hline
$112\%$ & $26.34$ & $26.2941$ & $1.742597\times 10^{-03}$ & $24.95$ & $25.3308$ & $1.526253\times 10^{-02}$ \\
\hline
\end{tabular}
}
\caption{EURO STOXX 50. {\tt DSabr\_I} model: $\sigma_{market}$ vs. $\sigma_{model}$.}
\label{tab:dynamicAllStoxx}
\end{table}

In Table \ref{tab:dynamicAllStoxxGPU_OpenMP}, the performance of the calibration procedure is illustrated. The speedup of the mono-GPU version is around $225$, while the $2$ GPUs version achieves a speedup up to $421$.

\begin{table}[!htb]
\centering
{\footnotesize

\begin{tabular}{|c|c|c|c|c|c|c|}
\hline
& CPU ($1$ thread) &$2$ threads &$4$ threads &$8$ threads & GPU & $2$ GPUs\\
\hline
Time ($s$) &$20162.40$ &  $10089.16$ & $5116.84$ & $2643.87$ &$89.95$ &  $47.81$\\
\hline
Speedup & - & $ 1.99$ & $3.94$ & $7.63$ & $224.15$ &$421.72$ \\
\hline
\end{tabular}
}
\caption{EURO STOXX 50. {\tt DSabr\_I} model: Performance of OpenMP vs. GPU versions, in single precision.}
\label{tab:dynamicAllStoxxGPU_OpenMP}
\end{table}

\item {\bf Calibration of the {\tt DSabr\_II} model using the technique {\tt T\_II}}

An asymptotic expression of implied volatility for the {\tt DSabr\_II} model is not available. So, for its calibration we use the technique {\tt T\_II} (see Section \ref{sec:calibracionGpu}). The calibration process is performed in prices. For the Monte Carlo pricing method we have considered $2^{20}$ paths with $\Delta t= 1/250$. In Table \ref{tab:dynamicSetCallibrateDataEURO STOXX 50_CII} the calibrated parameters are detailed. In this case, double precision computations have been used to ensure the convergence to the real minimum.

\begin{table}[!htb]
\centering
{\footnotesize
\begin{tabular}{|c|c|c|c|c|}
\hline
$\alpha=0.296790$ & $\beta=1.000000$ & $\rho_0=-0.360610$ & $\nu_0=0.000100$ & $a=15.0000000$\\
\hline
$b=15.000000$ & $d_{\rho}=-0.715716$ & $d_{\nu}=0.847244$ & $q_{\rho}=15.000000$ &$q_{\nu}=-8.969205$\\
\hline
\end{tabular}
}
\caption{EURO STOXX 50. {\tt DSabr\_II} model: Calibrated parameters.}
\label{tab:dynamicSetCallibrateDataEURO STOXX 50_CII}
\end{table}

In Figure \ref{fig:dynamicAllStoxx_CII}, the whole prices surface at maturities $3$, $6$, $12$ and $24$ months is shown. In Table \ref{tab:dynamicAllStoxx_CII} the market and model prices are compared. The mean relative error is $1.741038\times 10^{-2}$ and the maximum relative error is $5.344231\times 10^{-2}$.

\begin{figure}[h!]
\centering
\subfigure{\includegraphics[height=5cm]{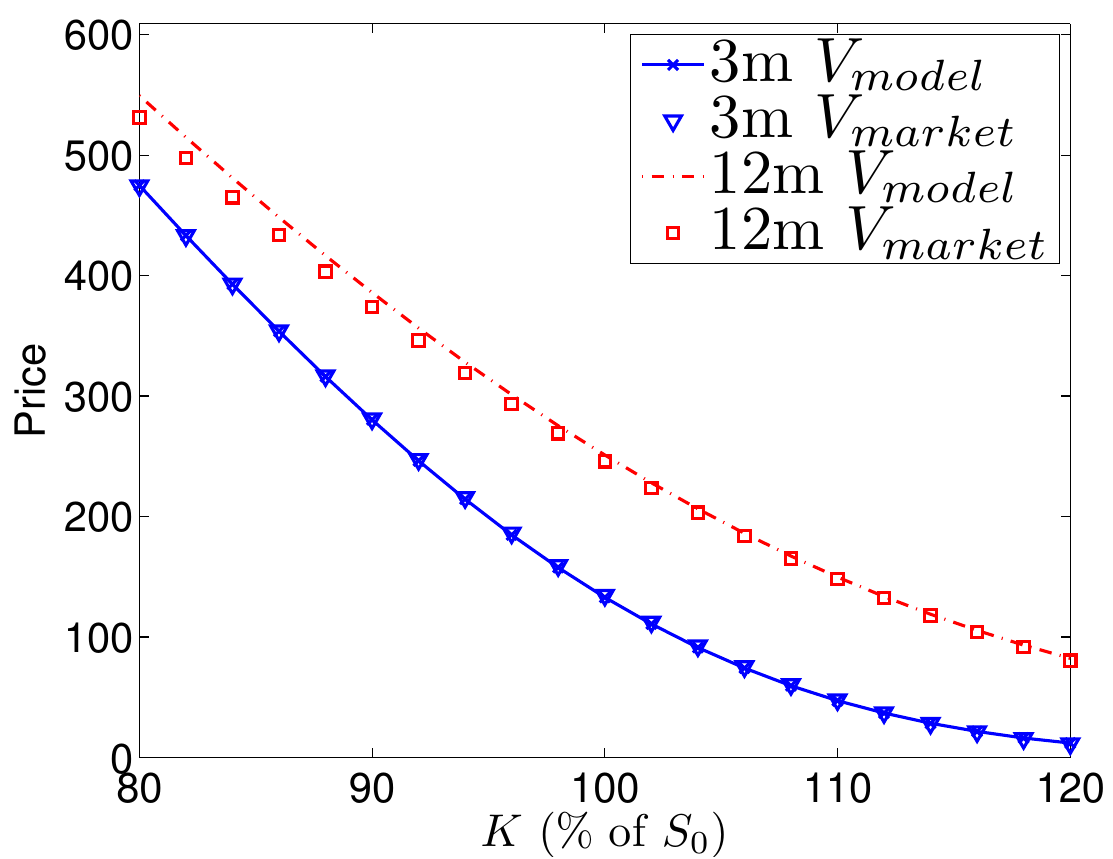}}
\subfigure{\includegraphics[height=5cm]{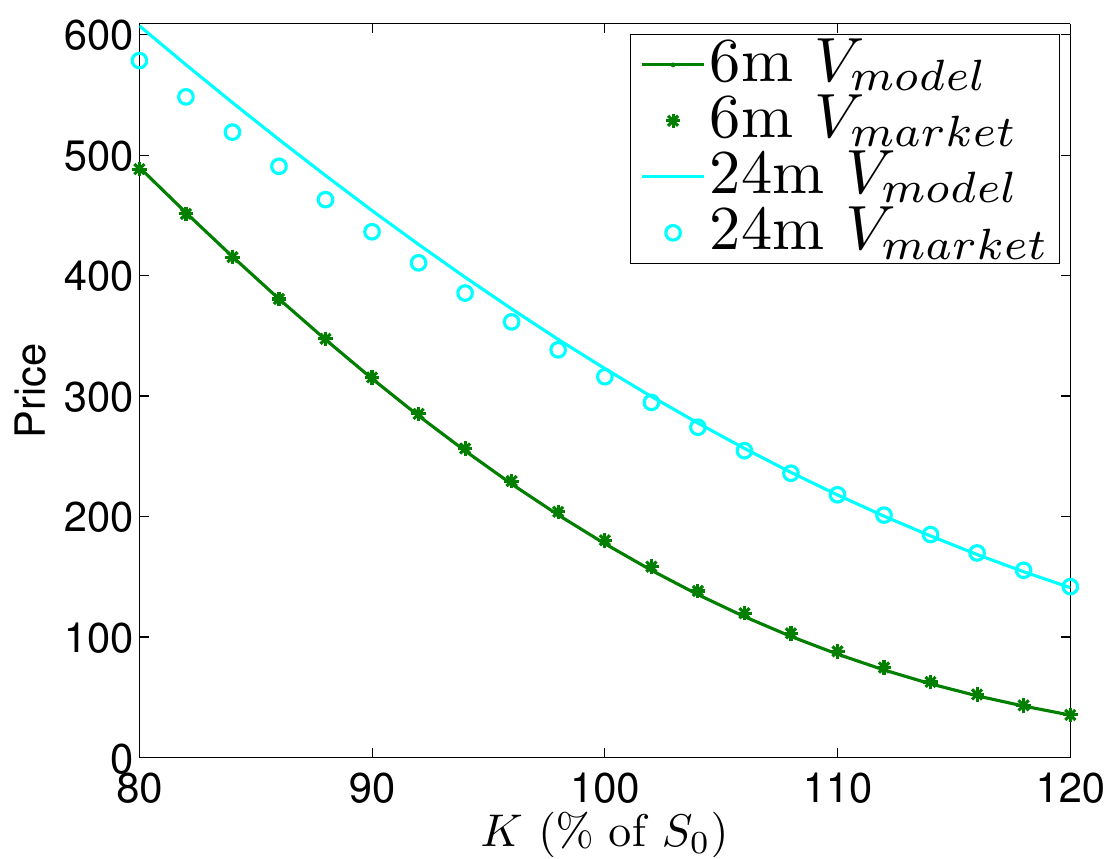}}
\caption{EURO STOXX 50. {\tt DSabr\_II} model: $V_{model}$ vs. $V_{market}$ for the whole prices surface. Maturities: 3 and 12 months (left), 6 and 24 months (right).}
\label{fig:dynamicAllStoxx_CII}
\end{figure}

\begin{table}[!htb]
\centering
{\footnotesize
\begin{tabular}{|c|| c|c|c ||c|c|c |}
 \hline
$K$ (\% of $S_0$) & \multicolumn{3}{|c|}{$3$ months }& \multicolumn{3}{|c|}{$6$ months } \\
\hline
\, & $V_{market}$ & $V_{model}$ & $\frac{|V_{market}-V_{model}|}{V_{market}}$
         & $V_{market}$ & $V_{model}$ &$\frac{|V_{market}-V_{model}|}{V_{market}}$ \\
\hline
$88\%$ & $316.679$ & $316.119$ & $1.769641\times 10^{-03}$ & $347.371$ & $346.830$ & $1.556418\times 10^{-03}$ \\
\hline
$100\%$ & $134.605$ & $132.952$ & $1.227935\times 10^{-02}$ & $180.353$ & $177.429$ & $1.621589\times 10^{-02}$ \\
\hline
$112\%$ & $37.252$ & $36.895$ & $9.572343\times 10^{-03}$ & $74.680$ & $72.682$ & $2.676232\times 10^{-02}$ \\
\hline
\hline
$K$ (\% of $S_0$) & \multicolumn{3}{|c|}{$12$ months } & \multicolumn{3}{|c|}{$24$ months}\\
\hline
\,& $V_{market}$ & $V_{model}$ & $\frac{|V_{market}-V_{model}|}{V_{market}}$
       & $V_{market}$ & $V_{model}$ & $\frac{|V_{market}-V_{model}|}{V_{market}}$\\
\hline
$88\%$ & $403.205$ & $416.516$ & $3.301443\times 10^{-02}$ & $463.037$ & $482.939$ & $4.298152\times 10^{-02}$ \\
\hline
$100\%$ & $245.905$ & $251.015$ & $2.077787\times 10^{-02}$ & $316.081$ & $322.821$ & $2.132105\times 10^{-02}$ \\
\hline
$112\%$ & $132.454$ & $133.687$ & $9.304021\times 10^{-03}$ & $201.189$ & $200.374$ & $4.048420\times 10^{-03}$ \\
\hline
\end{tabular}
}
\caption{EURO STOXX 50. {\tt DSabr\_II} model: $V_{market}$ vs. $V_{model}$.}
\label{tab:dynamicAllStoxx_CII}
\end{table}


In double precision, the computational time with one GPU is $37709.57$ seconds and with $2$ GPUs is $19520.73$ seconds, getting
a speedup around $1.93$. In this case the CPU computation cost results prohibitive. As expected, the computational time of the calibration process with the Monte Carlo pricing method is much higher than with the previous calibration procedures.


\end{enumerate}

\underline{{\tt DSabr\_I} vs. {\tt DSabr\_II}}\\
In order to compare the accuracy of {\tt DSabr\_II} and {\tt DSabr\_I} models, we compute the mean relative error of the {\tt DSabr\_I} model in prices, with Black-Scholes formula applied to $\sigma_{model}$. This error is $2.154846\times 10^{-2}$, so that {\tt DSabr\_II} model captures better the market dynamics.

\subsubsection{EUR/USD exchange rate}


In this section we consider the EUR/USD exchange rate, that expresses the amount of American Dollars equivalent to one Euro. In Tables \ref{eurusd} and \ref{eurusd_vola} of the Appendix \ref{Append:MarketData} the interests rates, dividend yields and volatility smiles for $3$, $6$, $12$ and $24$ months maturities are shown. The EUR/USD \textit{spot} rate is $S_0 = 1.2939 $ US dollars quoted in December of 2011. From now on we denote the EUR/USD foreign rate as EURUSD. Next, we present the results of calibrating the introduced models to these data.


\begin{enumerate}
 \item {\bf Calibration of the {\tt SSabr} model using the technique {\tt T\_I} }\\
By using the technique {\tt T\_I} and the asymptotic expression (\ref{vola_SABR}), the individual calibration of the {\tt SSabr} model has been carried out. The parameters in Table \ref{tab:staticCallibrateDataEURUSD_Individual} have been obtained. For them, in Figure \ref{fig:staticIndividualUsd} the market and the model volatilities are shown.

\begin{table}[!htb]
\centering
{\footnotesize
\begin{tabular}{|c|c|c | c |c|}
\hline
\, & $3$ months & $6$ months & $12$ months & $24$ months\\
\hline\hline
 $\alpha$ & $0.146859$ & $0.152825$ & $0.158210$ & $0.154572$\\
\hline
$\beta$ & $1.0$& $0.990518$ & $0.945088$ & $0.999993$\\
\hline
$\nu$ & $0.911966$ & $0.675457$ & $0.491647$ & $0.328907$\\
\hline
$\rho$&$-0.447718$& $-0.490521$ & $-0.511180$ & $-0.560022$ \\
\hline
\end{tabular}
}
\caption{EURUSD. {\tt SSabr} model: Calibrated parameters for each maturity.}
\label{tab:staticCallibrateDataEURUSD_Individual}
\end{table}

\begin{figure}[h!]
 \centering
 \subfigure{\includegraphics[height=5cm]{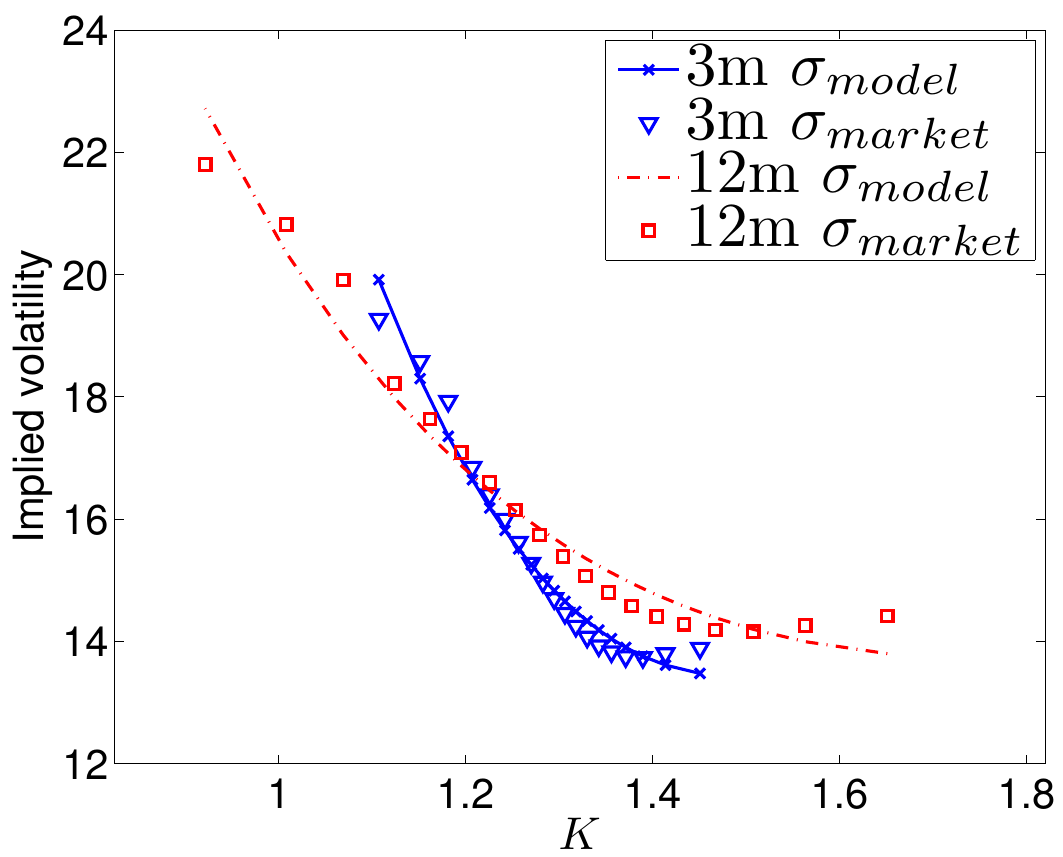}}
 \subfigure{\includegraphics[height=5cm]{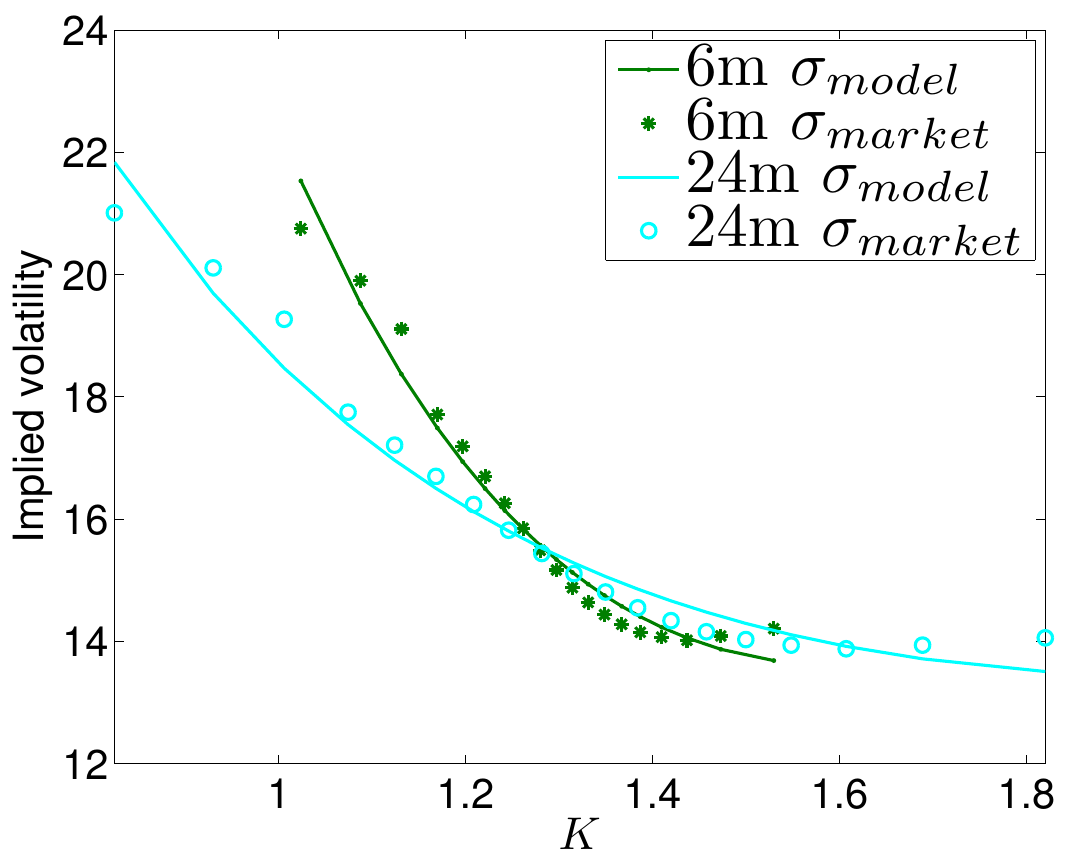}}
 \caption{EURUSD. {\tt SSabr} model: $\sigma_{model}$ vs. $\sigma_{market}$ for the whole volatility surface. Maturities: 3 and 12 months (left), 6 and 24 months (right).}
 \label{fig:staticIndividualUsd}
\end{figure}

In Table \ref{tab:static_24_EURUSD}, for $T=24$ months, the GPU performance is compared to one CPU. Time is measured in seconds and computation has been carried out in single precision. The speedup with $8$ OpenMP threads is near $8$, specifically $7.66$. As in EURO STOXX 50 calibration, 1 GPU speedup is around $190$ and $2$ GPUs around $333$.

\begin{table}[!htb]
\centering
{\footnotesize
\begin{tabular}{|c|c|c|c|c|c|c|}
\hline
& CPU ($1$ thread) & $2$ threads & $4$ threads & $8$ threads &GPU &$2$ GPUs\\
\hline
Time ($s$) &$4198.03$ & $2109.56$ & $1062.79$ & $548.04$ & $21.84$ & $12.58$\\
\hline
Speedup & - & $ 1.99$ & $3.95$ & $7.66$ & $192.19$ &$333.52$\\
\hline
\end{tabular}
}
\caption{EURUSD. {\tt SSabr} model: Performance of OpenMP vs. GPU versions, in single precision for $T = 24$ months.}
\label{tab:static_24_EURUSD}
\end{table}

\item {\bf Calibration of the {\tt DSabr\_I} model using the technique {\tt T\_I}}

{\tt DSabr\_I} model calibration to all strikes and maturities has been performed with the SA GPU version and formula \eqref{vola_SABRd}, with $\rho$ and $\nu$ given by (\ref{el_fun_1}).  In Table \ref{tab:dynamicSetCallibrateDataEURUSD} the calibrated parameters are detailed. In Figure \ref{fig:dynamicAllUsd}, the whole volatility surface at maturities $3$, $6$, $12$ and $24$ months is shown. Note that the dynamic SABR model captures correctly the volatility skew. In Table \ref{tab:dynamicAllEURUSD}, the market volatilities vs. the model ones \eqref{vola_SABRd} are shown. The mean relative error is $ 2.441714 \times 10^{-2}$ and the maximum relative error is $6.954307\times 10^{-2}$. In Table \ref{tab:dynamicAllUsdOpenMP_GPU}, the computational times and the speedups in single precision are shown. Note that for 1 GPU the speedup is around $240$, while for 2 GPUs is nearly $451$.

\begin{table}[!htb]
\centering
{\small
\begin{tabular}{|c|c|c|c|c|c|}
\hline
$\alpha=0.155464$ & $\beta=0.971908$ & $\rho_0=-0.642617$ & $\nu_0=0.800275$ & $a=0.001$ & $b=2.6093$\\
\hline
\end{tabular}
}
\caption{EURUSD. {\tt DSabr\_I} model: Calibrated parameters.}
\label{tab:dynamicSetCallibrateDataEURUSD}
\end{table}

\begin{figure}[h!]
\centering
    \subfigure{\includegraphics[height=5cm]{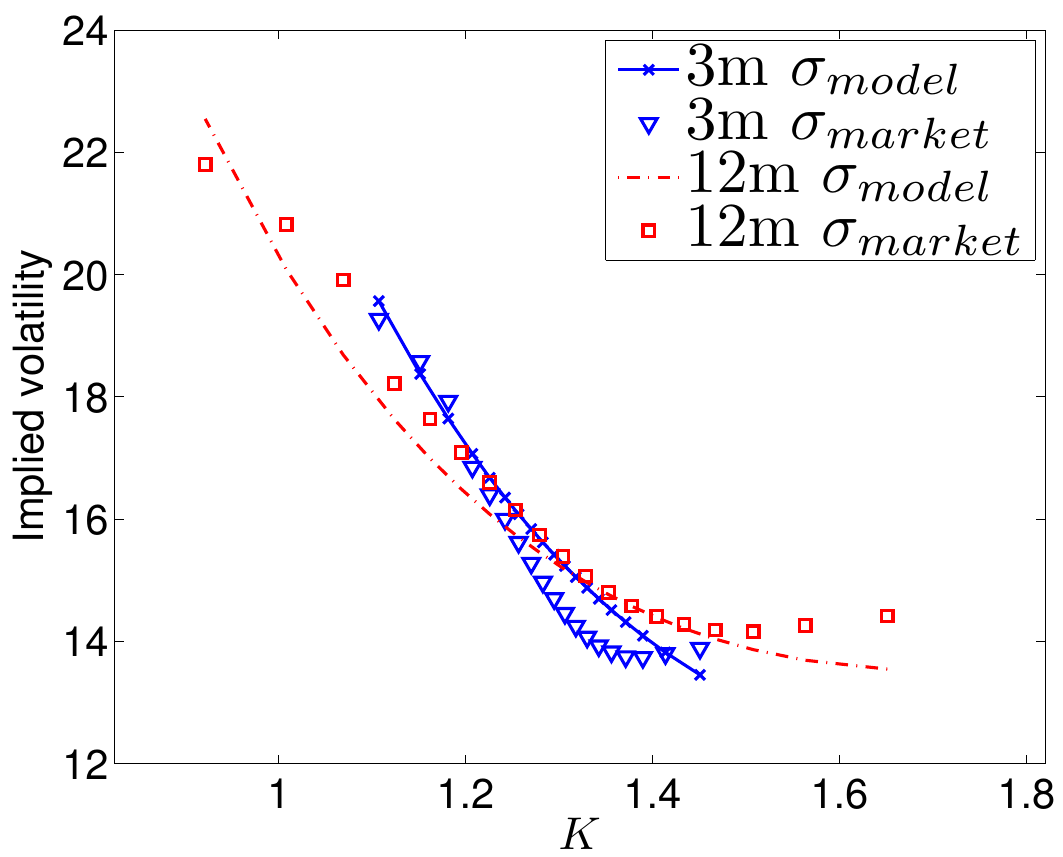}}
    \subfigure{\includegraphics[height=5cm]{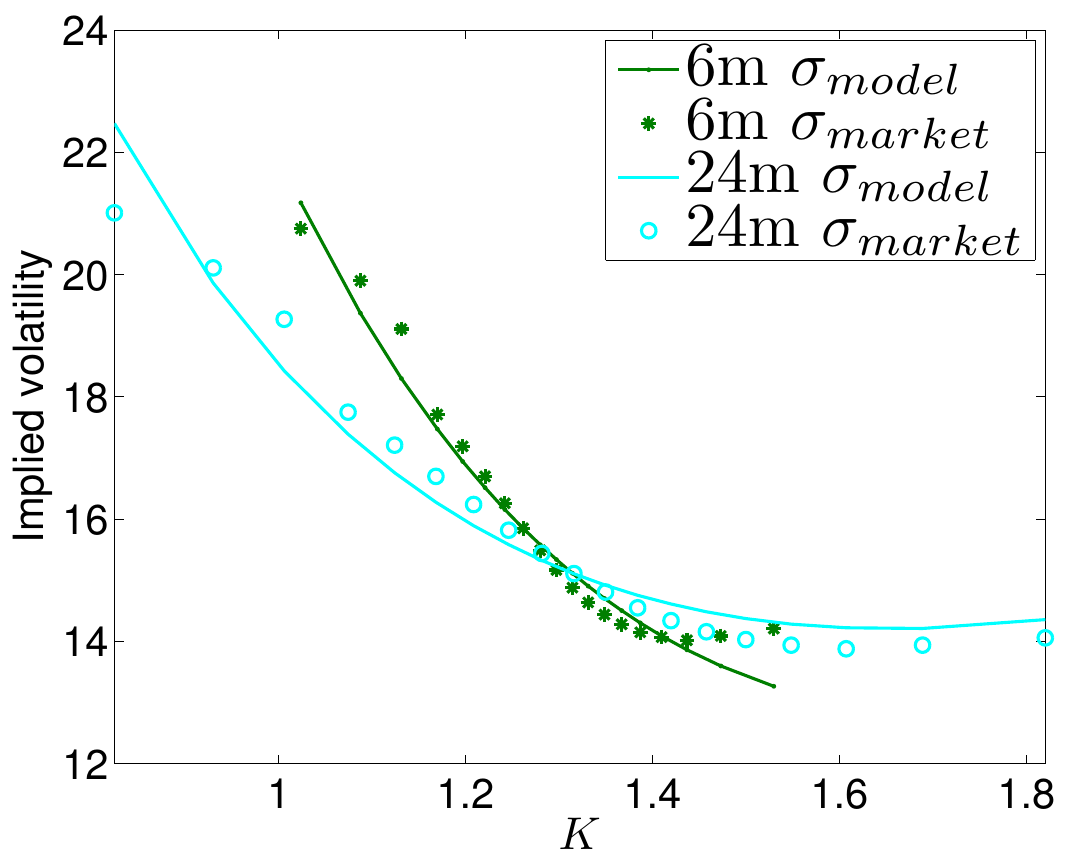}}
    \caption{EURUSD. {\tt DSabr\_I} model: $\sigma_{model}$ vs. $\sigma_{market}$ for the whole volatility surface. Maturities: 3 and 12 months (left), 6 and 24 months (right).}
  \label{fig:dynamicAllUsd}
\end{figure}

\begin{table}[!htb]
\centering
{\footnotesize
\begin{tabular}{|c| c|c|c ||c|c|c|c |}
 \hline
 \multicolumn{4}{|c|}{$3$ months }& \multicolumn{4}{|c|}{$6$ months } \\
\hline
$K$&  $\sigma_{market}$ & $\sigma_{model}$ & $\frac{|\sigma_{market}-\sigma_{model}|}{\sigma_{market}}$
        &$K$ & $\sigma_{market}$ & $\sigma_{model}$ &$\frac{|\sigma_{market}-\sigma_{model}|}{\sigma_{market}}$ \\
\hline
$1.2075$ & $16.85$ & $17.0683$ & $1.295549\times 10^{-02}$ & $1.1700$ & $17.71$ & $17.4751$ & $1.326369\times 10^{-02}$ \\
\hline
$1.2950$ & $14.70$ & $15.4197$ & $4.895918\times 10^{-02}$ &$1.2975$ & $15.17$ & $15.3398$ & $1.119314\times 10^{-02}$ \\
\hline
$1.3715$ & $13.75$ & $14.3171$ & $4.124364\times 10^{-02}$ &  $1.4099$ & $14.07$ & $14.0914$ & $1.520967\times 10^{-03}$ \\
\hline
\hline
 \multicolumn{4}{|c|}{$12$ months } & \multicolumn{4}{|c|}{$24$ months}\\
\hline
$K$& $\sigma_{market}$ & $\sigma_{model}$ & $\frac{|\sigma_{market}-\sigma_{model}|}{\sigma_{market}}$
      &$K$ & $\sigma_{market}$ & $\sigma_{model}$ & $\frac{|\sigma_{market}-\sigma_{model}|}{\sigma_{market}}$\\
\hline
$1.1240$ & $18.22$ & $17.6324$ & $3.225027\times 10^{-02}$ & $1.0746$ & $17.75$ & $17.3887$ & $2.035493\times 10^{-02}$ \\
\hline
$1.3043$ & $15.39$ & $15.2020$ & $1.221572\times 10^{-02}$ &  $1.3161$ & $15.11$ & $15.1075$ & $1.654533\times 10^{-04}$ \\
\hline
$1.4673$ & $14.19$ & $14.0396$ & $1.059901\times 10^{-02}$ &  $1.5485$ & $13.94$ & $14.2853$ & $2.477044\times 10^{-02}$ \\
\hline
\end{tabular}
}
\caption{EURUSD. {\tt DSabr\_I} model: $\sigma_{market}$ vs. $\sigma_{model}$.}
\label{tab:dynamicAllEURUSD}
\end{table}

\begin{table}[!htb]
\centering
{\footnotesize
\begin{tabular}{|c|c|c|c|c|c|c|}
\hline
&CPU ($1$ thread) & $2$ threads & $4$ threads &$8$ threads &GPU &$2$ GPUs\\
\hline
Time ($s$) &$16793.41$ &$8389.14 $&$4240.11 $& $2204.20 $ & $69.73$ & $37.16$ \\
\hline
Speedup  & - & $2.00 $& $3.96 $ & $ 7.62 $ & $240.83$ & $451.92$\\
\hline
\end{tabular}
}
\caption{EURUSD. {\tt DSabr\_I} model: Performance of OpenMP vs. GPU versions, in single precision.}
\label{tab:dynamicAllUsdOpenMP_GPU}
\end{table}

\item {\bf Calibration of the {\tt DSabr\_II} model using technique {\tt T\_II}}

Analogously to the previous Section \ref{subsub_eurostoxx}, an asymptotic expression for implied volatility in the {\tt DSabr\_II} model is not available. The calibration has been carried with Technique II (see  Section \ref{sec:calibracionGpu}). In this case, calibration is performed in prices and using double precision. The set of calibrated parameters is detailed in Table \ref{tab:dynamicSetCallibrateDataEURUSD_CII}.
\begin{table}[!htb]
\centering
{\small
\begin{tabular}{|c|c|c|c|c|}
\hline
$\alpha=0.154037$ & $\beta=1.000000$ & $\rho_0=-0.693682$ & $\nu_0=7.541424$ & $a=0.000000$\\
\hline
 $b=150.000000$  &$d_{\rho}=-0.200342$ &  $d_{\nu}=0.339807$ & $q_{\rho}=0.345973$ &$q_{\nu}=-0.992551$\\
\hline
\end{tabular}
}
\caption{EURUSD. {\tt DSabr\_II} model: Calibrated parameters.}
\label{tab:dynamicSetCallibrateDataEURUSD_CII}
\end{table}

Figure \ref{fig:dynamicAllUsd_CII} shows the comparison between market and model prices. The maximum relative error is $1.418863\times 10^{-1}$ and the mean
relative error is $2.192849 \times 10^{-2}$. In Table \ref{tab:dynamicAllUsd_CII} the market and model prices for some strikes are shown.

\begin{figure}[h!]
\centering
\subfigure{\includegraphics[height=5cm]{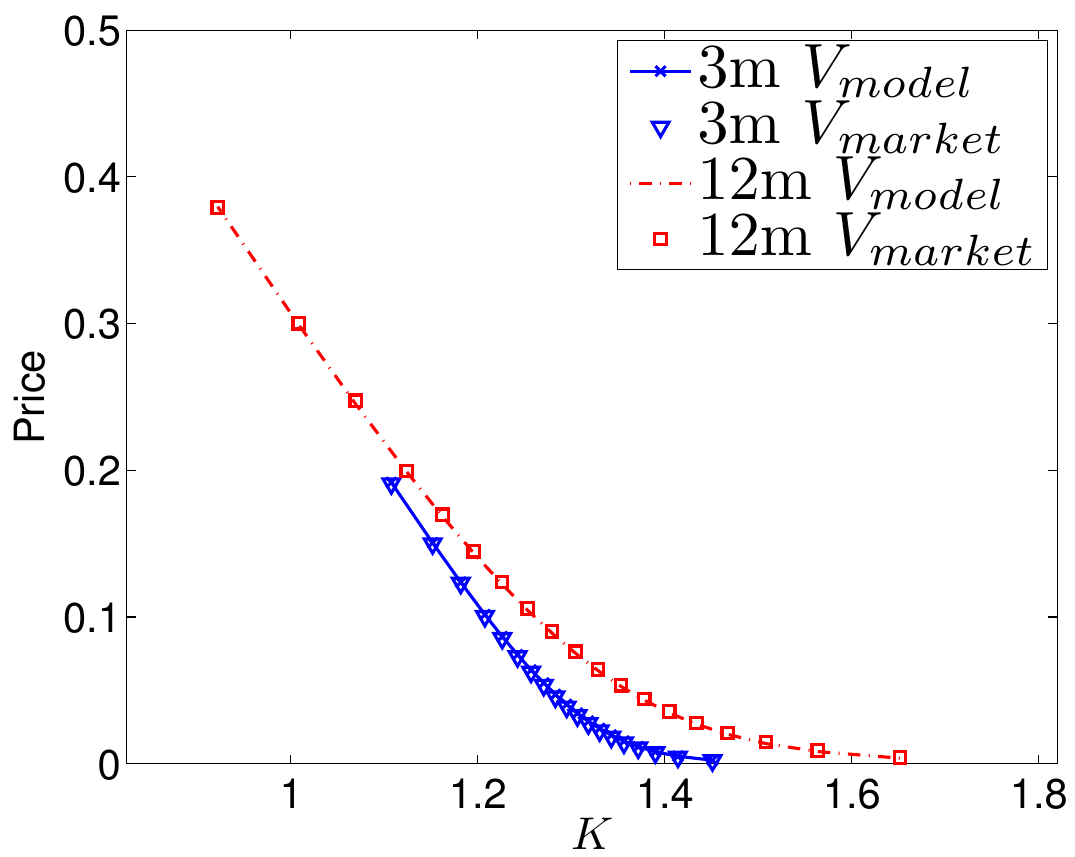}}
\subfigure{\includegraphics[height=5cm]{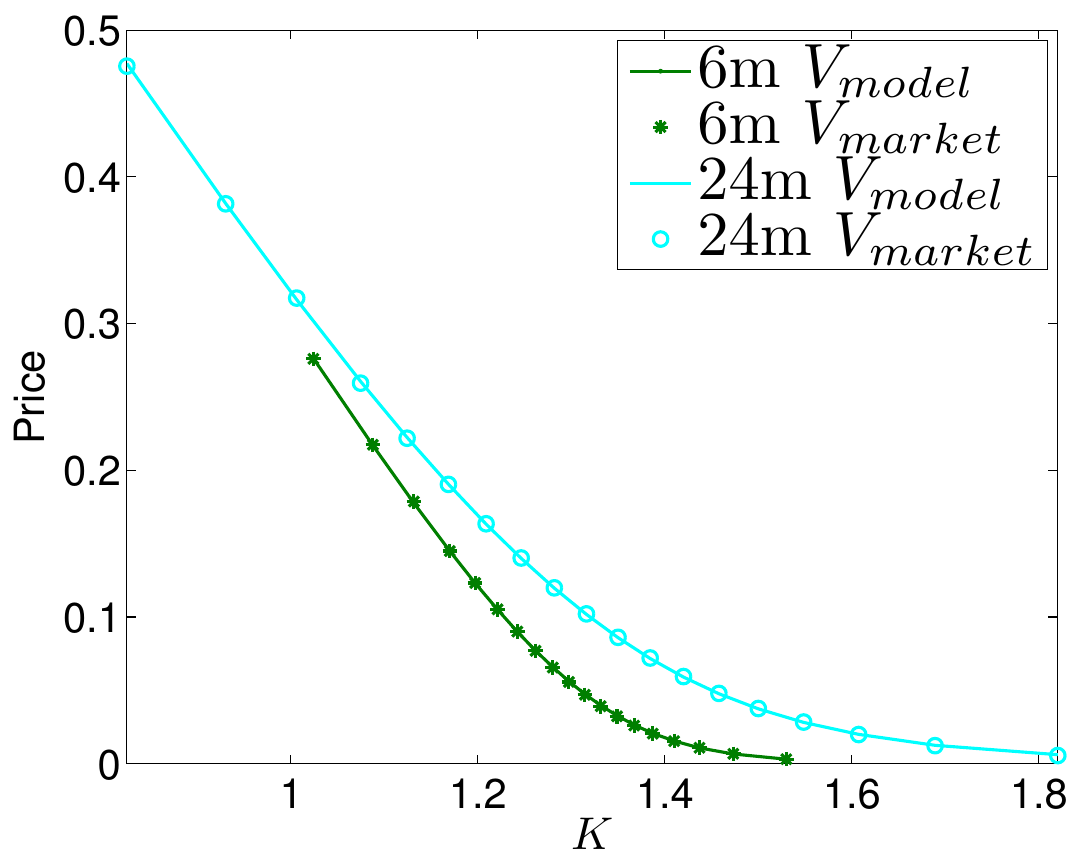}}
\caption{EURUSD. {\tt DSabr\_II} model: $V_{model}$ vs. $V_{market}$ for the whole prices surface. Maturities: 3 and 12 months (left), 6 and 24 months (right).}
\label{fig:dynamicAllUsd_CII}
\end{figure}

\begin{table}[!htb]
\centering
{\small
\begin{tabular}{|c| c|c|c ||c|c|c|c |}
 \hline
 \multicolumn{4}{|c|}{$3$ months }& \multicolumn{4}{|c|}{$6$ months } \\
\hline
$K$&  $V_{market}$ & $V_{model}$ & $\frac{|V_{market}-V_{model}|}{V_{market}}$
        &$K$ & $V_{market}$ & $V_{model}$ &$\frac{|V_{market}-V_{model}|}{V_{market}}$ \\
\hline
$1.2075$ & $0.100489$ & $0.101712$ & $1.217204\times 10^{-02}$ & $1.1700$ & $0.144905$ & $0.144950$ & $3.072116\times 10^{-04}$ \\
\hline
$1.2950$ & $0.038794$ & $0.040476$ & $4.335060\times 10^{-02}$ & $1.2975$ & $0.055770$ & $0.056139$ & $6.611759\times 10^{-03}$ \\
\hline
$1.3715$ & $0.010770$ & $0.011697$ & $8.603287\times 10^{-02}$ & $1.4099$ & $0.015472$ & $0.015539$ & $4.271299\times 10^{-03}$ \\
\hline
\hline
 \multicolumn{4}{|c|}{$12$ months } & \multicolumn{4}{|c|}{$24$ months}\\
\hline
$K$& $V_{market}$ & $V_{model}$ & $\frac{|V_{market}-V_{model}|}{V_{market}}$
      &$K$ & $V_{market}$ & $V_{model}$ & $\frac{|V_{market}-V_{model}|}{V_{market}}$\\
\hline
$1.1240$ & $0.199490$ & $0.198897$ & $2.971219\times 10^{-03}$ & $1.0746$ & $0.259398$ & $0.260539$ & $4.399383\times 10^{-03}$ \\
\hline
$1.3043$ & $0.076379$ & $0.075409$ & $1.269954\times 10^{-02}$ & $1.3161$ & $0.102106$ & $0.101766$ & $3.330301\times 10^{-03}$ \\
\hline
$1.4673$ & $0.020951$ & $0.020010$ & $4.495205\times 10^{-02}$ & $1.5485$ & $0.028409$ & $0.028189$ & $7.756103\times 10^{-03}$ \\
\hline
\end{tabular}
}
\caption{EURUSD. {\tt DSabr\_II} model: $V_{market}$ vs. $V_{model}$. }
\label{tab:dynamicAllUsd_CII}
\end{table}


The performance of calibration using 1 GPU is $35033.72$ seconds and with 2 GPUs is
equal to $19143.35$ seconds; getting a speedup around $1.83$ times. We do not make a comparison between CPU and GPU because the computation time in CPU results to be prohibitive.

%

\end{enumerate}

\underline{{\tt DSabr\_I} vs. {\tt DSabr\_II}}\\
The accuracy of {\tt DSabr\_II} model with respect to {\tt DSabr\_I} model is compared using the mean relative error in prices. The {\tt DSabr\_I} model error is $3.647441\times 10^{-2}$. Therefore, again the {\tt DSabr\_II} model captures better the market dynamics.

\subsubsection{Calibration test: pricing European options}
In order to validate the correct calibration of model parameters, we price European options with
the {\tt DSabr\_I} model and Monte Carlo pricing method. The price is denoted with $V_{\sigma_{model}}$. Note that for the {\tt DSabr\_II} model this task is redundant, since calibration is also carried out with the same Monte Carlo method.

\begin{figure}[!htb]
  \begin{center}
    \subfigure{\includegraphics[height=5cm]{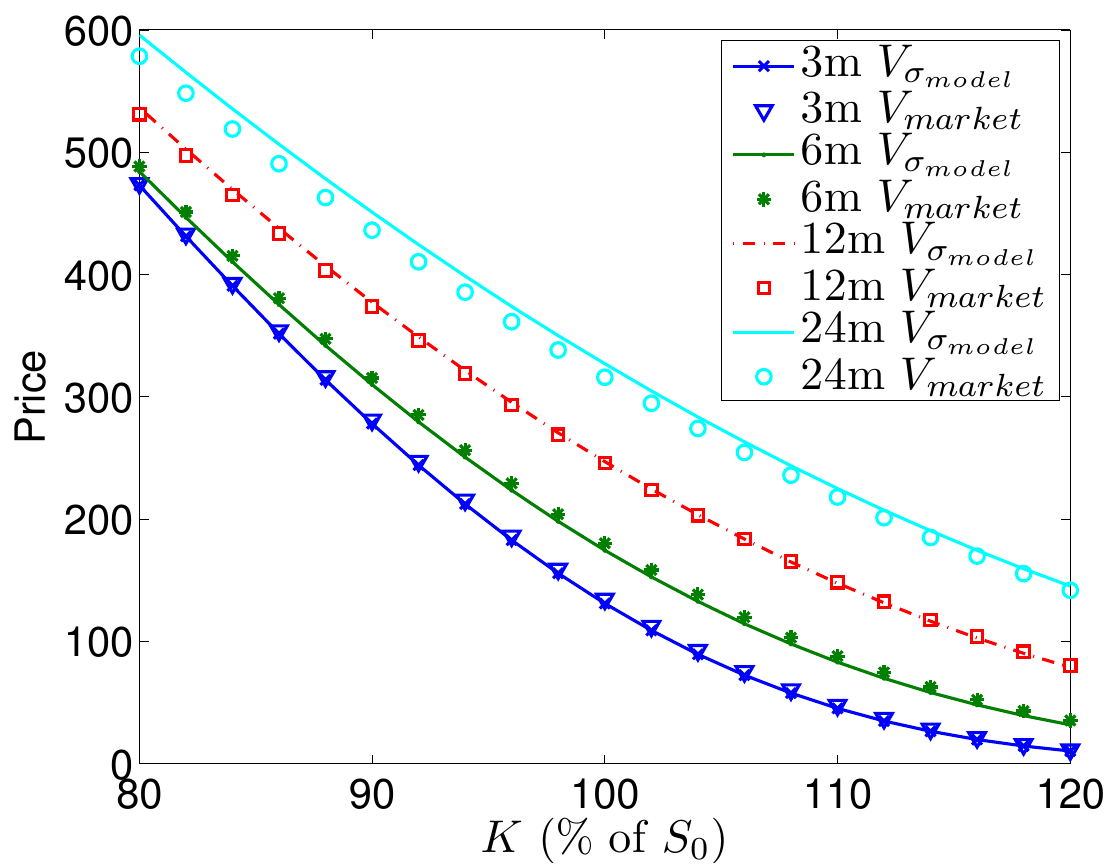}}
    \subfigure{\includegraphics[height=5cm]{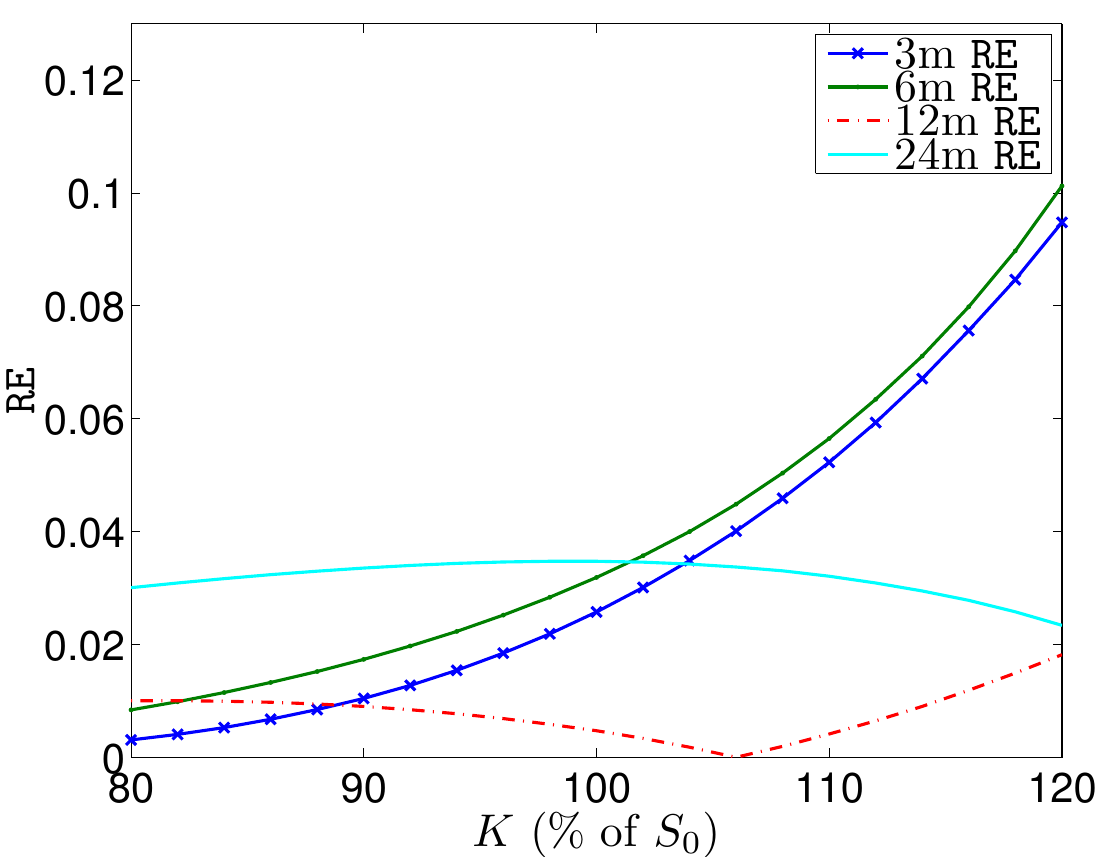}}
    \caption{EURO STOXX 50. {\tt DSabr\_I} model for pricing European options. Prices (left) and relative errors (right).}
    \label{fig:pricingDynamicAll_x}
  \end{center}
\end{figure}

\noindent\textbf{EURO STOXX 50:} In Figure \ref{fig:pricingDynamicAll_x}, the comparison between market prices (calculated with Black-Scholes formula and market volatilities) and model prices (calculated with the expression \eqref{vola_SABRd} and parameters of Table \ref{tab:dynamicSetCallibrateDataEURO STOXX 50} in the Black-Scholes formula) are plotted. The mean relative error is $2.840228 \times 10^{-2}$ and the
maximum relative error is $1.008734 \times 10^{-1}$. Pricing all options with GPU takes $0.257248$ seconds.

\begin{figure}[!htb]
  \begin{center}
    \subfigure{\includegraphics[height=5cm]{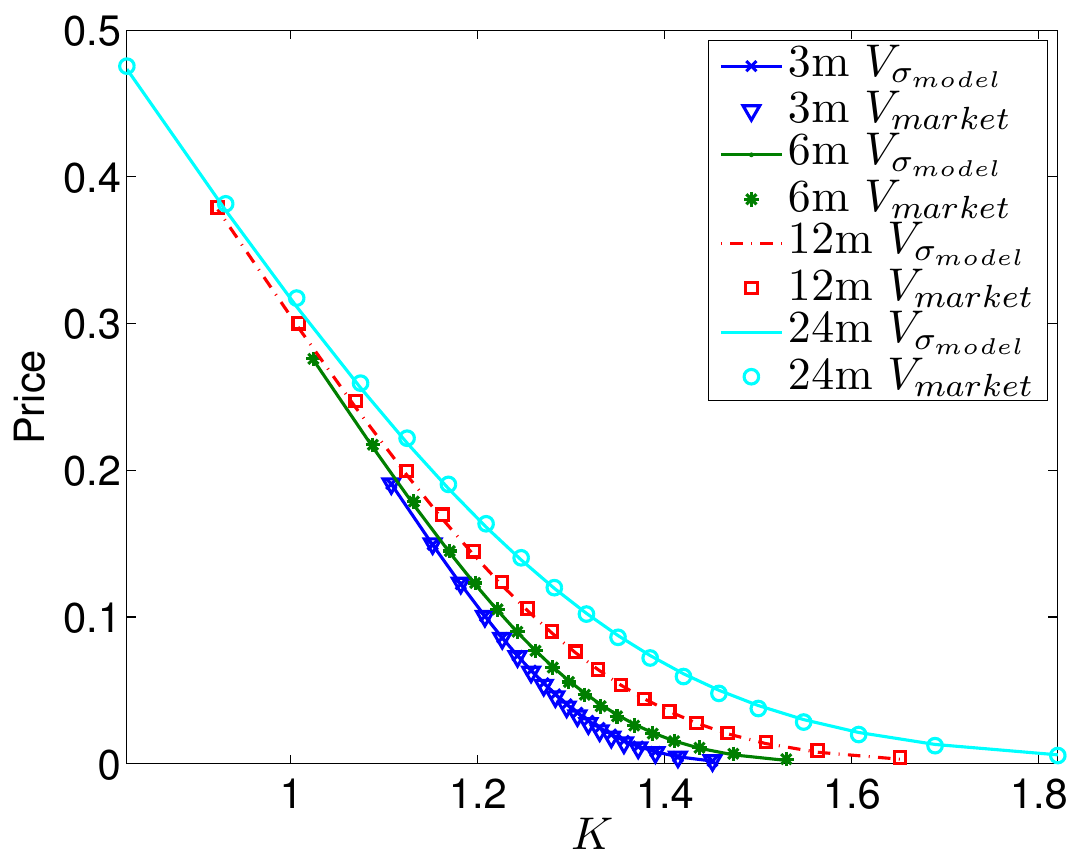}}
     \subfigure{\includegraphics[height=5cm]{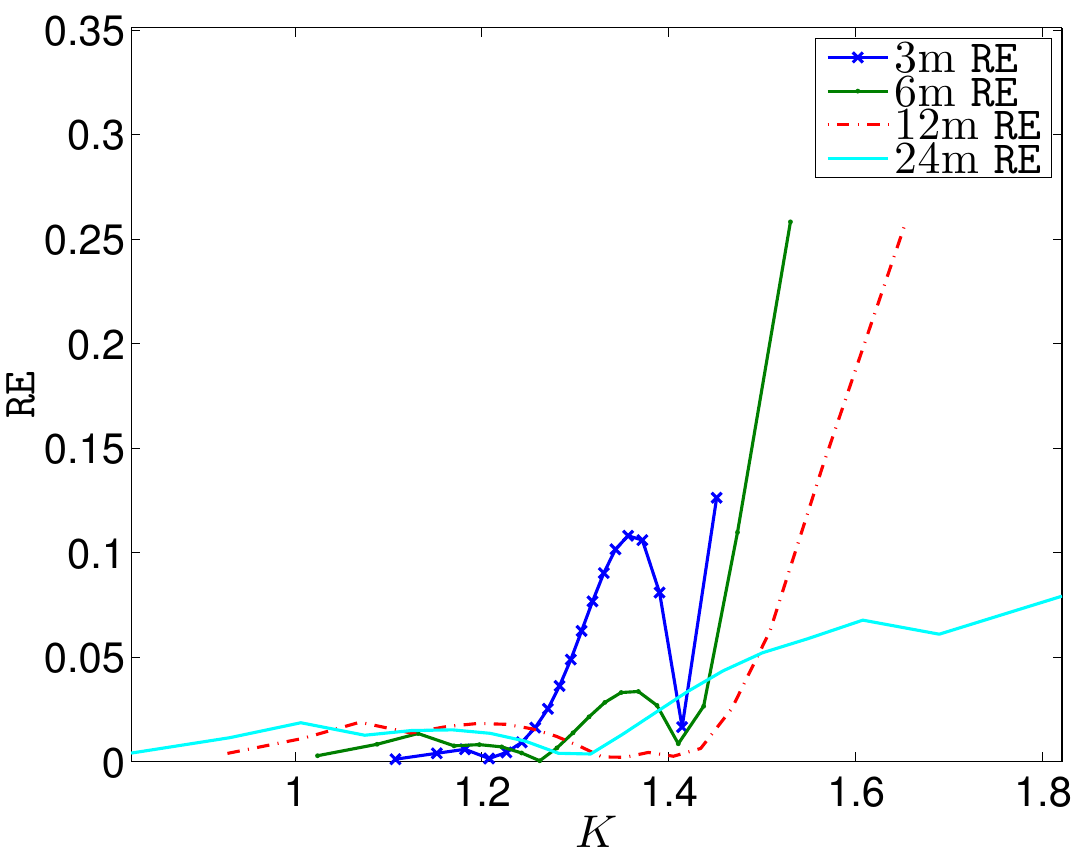}}
    \caption{EURUSD. {\tt DSabr\_I} model for pricing European options. Prices (left) and relative errors (right).}
    \label{fig:pricingDynamicAll}
  \end{center}
\end{figure}

 \noindent\textbf{EURUSD:} Analogously of the EURO STOXX 50 case, in Figure \ref{fig:pricingDynamicAll}, the market prices vs. the model ones are shown, using the parameters of Table \ref{tab:dynamicSetCallibrateDataEURUSD}. Furthermore, relative errors are shown. The mean relative error is $3.577857 \times 10^{-02}$ and the maximum relative error is $2.582023 \times 10^{-01}$. Pricing all options in GPU takes $0.248013$ seconds.


\subsection{Pricing a cliquet option}

In this section the pricing of a cliquet option on the EUR/USD exchange rate is described. For this purpose, we use the GPU Monte Carlo method and consider {\tt DSabr\_I} and {\tt DSabr\_II} models. For the first one, we choose the
parameters of Table \ref{tab:dynamicSetCallibrateDataEURUSD} and for the second one those in Table \ref{tab:dynamicSetCallibrateDataEURUSD_CII}.

Cliquet options are options where the strike price is periodically reset several times before final expiration date \cite{cliquetWilmot2002}. At the resetting date, the option will expire worthless if the current security price is below the strike price, and the strike price will be reset to this lower price. If the security price at resetting date is higher than the strike, the investor will earn the difference and the strike price will be reset to this higher price. Thus, a cliquet option is equivalent to a series of forward-starting at-the-money options with local limits, cap and floor. The option payoff function is:

\begin{equation}				
Cliquet\big(F_{l},C_{l},\{d_{i}\}_{i=1\ldots D},\{S_{d_{i}}\}_{i=1\ldots D}\big)=
\sum_{i=2}^{D}\max\left(\min\left(\frac{S_{d_{i}}-S_{d_{i-1}}}{S_{d_{i-1}}},C_{l}\right),F_{l}\right),
\label{eq_cliquet}
\end{equation}
\noindent where $C_l$ and $F_l$ are the local cap and floor limits, respectively, $d_{i}$ denotes a resetting date and $S_{d_{i}}$ is the underlying price at time $d_{i}$. Thus, $R_i=\frac{S_{d_{i}}-S_{d_{i-1}}}{S_{d_{i-1}}}$ is known as the \textit{return} between dates $d_{i-1}$ and $d_i$.

Usually, global limits can also be added, so that the payoff function is:
\begin{equation}
 Cliquet_{global}(C_{g},F_{g},Cliquet)=\max\big(C_{g},\min(F_{g},Cliquet)\big),
\label{eq:cliquet_global}
\end{equation}
\noindent where $C_g$ is global cap, $F_g$ global floor and $Cliquet$ is given by \eqref{eq_cliquet}.


 We consider the following data: $T=12$ months, $D=4$, $d_i= \frac{i \times T}{D}$, $F_l=-0.02$, $C_l=0.02$, $F_g=0$, $C_g=0.2$ and the payoff \eqref{eq:cliquet_global}. We price the cliquet option with {\tt DSabr\_I} and {\tt DSabr\_II} models, in single
precision. By using the {\tt DSabr\_I} model and the calibrated parameters in
Table \ref{tab:dynamicSetCallibrateDataEURUSD}, the simulated
cliquet option price is \textbf{$V_{\sigma_{model}}$}$=0.098289$, while when using {\tt DSabr\_II} model and the parameters in Table \ref{tab:dynamicSetCallibrateDataEURUSD_CII} then the
simulated price is \textbf{$V_{MC}$}$=0.120846$.  Note that $|V_{\sigma_{model}}-V_{MC}|=0.022557$ and  the execution time is $0.495199$ seconds in both cases.



\section{Conclusions}
Static or dynamic SABR models should be chosen depending on the particular financial derivative and the available market data. In both cases, the calibration of parameters can be carried out either by using an asymptotic implied volatility formula or a Monte Carlo simulation method. When using standard hardware tools, due to the high computational time with Monte Carlo strategy, the formula is mainly used. Nevertheless, the recently increasing use of GPU technology for scientific computing can also be extended to the algorithms involved in the calibration procedure. In the present paper we propose the application of this technology to the SA global optimization algorithm  and to the Monte Carlo simulation for the calibration in static and dynamic SABR models. In this GPU setting, the calibration by Monte Carlo is affordable in terms of computational time and the cost of calibration with the asymptotic formula can be highly reduced. In order to illustrate the performance of our GPU implementations, the calibration of static and dynamic SABR models for EURO STOXX 50 index and EUR/USD exchange rate have been carried out with asymptotic formula and Monte Carlo method. Once the parameters have been calibrated, a cliquet option on EUR/USD has been priced.  For the dynamic SABR model we propose an original more general expression for the functional parameters that reveals specially well suited for a EUR/USD exchange rate market data set.

Numerical results illustrate the expected behavior of both SABR models and the accuracy of the calibration. In terms of computational time, if we use the formula then the achieved speedup with respect to CPU computation is around $200$ with 1 GPU. Also, it is illustrated that GPU technology allows the use of Monte Carlo simulation for calibration purposes, the corresponding computational time with CPU being unaffordable. \\
\pagebreak
\appendix
{\LARGE \textbf{\appendixname}}

\section{Expression of implied volatility in the general case}\label{anexo_funciones_casogeneral}

For the \textbf{Case II} of Section \ref{elec_func} (general case), using Mathematica, the functions $\nu_1^2$, $\nu_2^2$, $\eta_1$ and $\eta_2^2$ given by (\ref{integrales_SABRd}) have the following expressions:
$$
\renewcommand\arraystretch{1.0}
{\small{
\begin{array}{ll}
\nu_1^2(T) = & \displaystyle \frac{1}{4 b^5 T^3} \Bigg[9q_{\nu}^2+6 b^4\nu_0 (4d_{\nu}+\nu_0) T^2+4 b^5d_{\nu}^2 T^3+9 bq_{\nu} (16d_{\nu}+\nu_0-q_{\nu} T)\\
& \displaystyle +  6 b^3 T (-\nu_0 (8d_{\nu}+\nu_0)+(4d_{\nu}+\nu_0)q_{\nu} T)+3 b^2 \left(\nu_0 (16d_{\nu}+\nu_0)-4 (8d_{\nu}+\nu_0)q_{\nu} T+q_{\nu}^2 T^2\right)\\
& \displaystyle - 3 e^{-2 b T} \bigg(3q_{\nu}^2+3 bq_{\nu} \Big(16d_{\nu} e^{b T}+\nu_0+q_{\nu} T\Big)+b^2 (\nu_0+q_{\nu} T) \Big(16d_{\nu} e^{b T}+\nu_0+q_{\nu} T\Big)\bigg) \Bigg],\\
\\
\nu_2^2(T) = &   \displaystyle \frac{1}{4 b^5 T^3}e^{-2 b T} \bigg\{18q_{\nu}^2+6 b^3 T (\nu_0+q_{\nu} T)^2+6 b^2 (\nu_0+q_{\nu} T) (\nu_0+3q_{\nu} T)\\
&  \displaystyle +9 bq_{\nu} (2\nu_0+3q_{\nu} T)+e^{2 b T} \Big[-18q_{\nu}^2+6 b^3\nu_0 (8d_{\nu}+\nu_0) T+4 b^5d_{\nu}^2 T^3   \\
& \displaystyle  +9 bq_{\nu} (-32d_{\nu}-2\nu_0+q_{\nu} T)+6 b^2 \big(-\nu_0 (16d_{\nu}+\nu_0)+2 (8d_{\nu}+\nu_0)q_{\nu} T\big)\Big]\\
& \displaystyle  +48 bd_{\nu} e^{b T} \Big(6q_{\nu}+b \big(\nu_0 (2+b T)+q_{\nu} T (4+b T)\big)\Big)\bigg\},
\end{array}
}}
$$
$$
\renewcommand\arraystretch{1.05}
{\small{
\begin{array}{ll}
\eta_1(T) =& \displaystyle \frac{2}{T^2} \Bigg\{-\frac{2d_{\nu}q_{\rho}}{a^3}-\frac{2d_{\rho}q_{\nu}}{b^3}-\frac{6q_{\rho}q_{\nu}}{(a+b)^4}+
         \frac{d_{\rho}\nu_0 T}{b}+\frac{d_{\nu}\rho_0 T}{a}+\frac{a^3\nu_0\rho_0 T}{(a+b)^4}+\frac{b^3\nu_0\rho_0 T}{(a+b)^4}   \\
 &  \displaystyle +\frac{d_{\nu} (-\rho_0+q_{\rho} T)}{a^2}  +\frac{d_{\rho} (-\nu_0+q_{\nu} T)}{b^2}+\frac{a^2 (-\nu_0\rho_0+\nu_0q_{\rho} T+q_{\nu}\rho_0 T)}{(a+b)^4} \\
& \displaystyle -\frac{2 a \big(\nu_0q_{\rho}+q_{\nu} (\rho_0-q_{\rho} T)\big)}{(a+b)^4}\\
& \displaystyle +\frac{b \big[-2\nu_0 (q_{\rho}+a\rho_0)+a\nu_0 (2q_{\rho}+3 a\rho_0) T+2q_{\nu} \big(q_{\rho} T+\rho_0 (-1+a T)\big)\big]}{(a+b)^4}\\
& \displaystyle   +\frac{b^2 \big[q_{\nu}\rho_0 T+\nu_0 \big(q_{\rho} T+\rho_0 (-1+3 a T)\big)\big]}{(a+b)^4}\\
& \displaystyle  +  \frac{1}{2 a^3 b^3 (a+b)^4} e^{-(a+b) T} \bigg\{4 b^7d_{\nu} e^{b T}q_{\rho}+8 a^2 b^5d_{\nu} e^{b T} \big(b\rho_0+q_{\rho} (3+b T)\big)   \\
& \displaystyle  +2 a b^6d_{\nu} e^{b T} \big( b\rho_0+q_{\rho} (8+b T) \big)+a^7d_{\rho} e^{a T} \left(4q_{\nu}+b^3d_{\nu} e^{b T} T^2+2 b (\nu_0+q_{\nu} T)\right)\\
&\displaystyle  +4 a^6 bd_{\rho} e^{a T} \left(4q_{\nu}+b^3d_{\nu} e^{b T} T^2+2 b (\nu_0+q_{\nu} T)\right)\\
&\displaystyle +2 a^5 b^2 \left[12d_{\rho} e^{a T}q_{\nu}+3 b^3d_{\rho}d_{\nu} e^{(a+b) T} T^2+6 bd_{\rho} e^{a T} (\nu_0+q_{\nu} T) \right.\\
& \displaystyle \left. +  b (\rho_0+q_{\rho} T) \left(d_{\nu} e^{b T}+\nu_0+q_{\nu} T\right)\right]+4 a^4 b^3 \left[d_{\nu} e^{b T}q_{\rho}+
             \nu_0q_{\rho}+4d_{\rho} e^{a T}q_{\nu}+q_{\nu}\rho_0+2q_{\rho}q_{\nu} T  \right. \\
&\displaystyle \left. +b^3d_{\rho}d_{\nu} e^{(a+b) T} T^2+2 bd_{\rho} e^{a T} (\nu_0+q_{\nu} T)+b (\rho_0+q_{\rho} T) \left(2d_{\nu} e^{b T}+\nu_0+q_{\nu} T\right)\right]\\
& \displaystyle +a^3 b^3 \bigg[12q_{\rho}q_{\nu}+b^4d_{\rho}d_{\nu} e^{(a+b) T} T^2+4 b \left(4d_{\nu} e^{b T}q_{\rho}+\nu_0q_{\rho}+
           q_{\nu} \left(d_{\rho} e^{a T}+\rho_0+2q_{\rho} T\right)\right)\\
& \displaystyle   +2 b^2 \Big[d_{\rho} e^{a T} (\nu_0+q_{\nu} T)+(\rho_0+q_{\rho} T) \left(6d_{\nu} e^{b T}+\nu_0+q_{\nu} T\right)
\Big]   \bigg]\bigg\}\Bigg\},
\end{array}
}}
$$

$$
\renewcommand\arraystretch{1.05}
{\small{
\begin{array}{ll}
\eta_2^2(T) =& \displaystyle \frac{12}{T^{4}}\displaystyle\int_{0}^{T}\int_0^t \Bigg\{\frac{1}{a^2 b^2 (a+b)^3}e^{-(a+b) s} \bigg\{ b^5d_{\nu} e^{b s} \left(-1+e^{a s}\right)q_{\rho} 6\\
& \displaystyle  -a b^4d_{\nu} e^{b s} \big(3q_{\rho}+b\rho_0-e^{a s} (3q_{\rho}+b\rho_0)+bq_{\rho} s\big)\\
& \displaystyle   +a^5d_{\rho} e^{a s} \left(-q_{\nu}-b (\nu_0+q_{\nu} s)+e^{b s} \big(q_{\nu}+b (\nu_0+bd_{\nu} s) \big)\right)+a^3 b^2 \left[-q_{\nu} \left(3d_{\rho} e^{a s}+\rho_0\right)   \right. \\
& \displaystyle   +e^{(a+b) s} \left((d_{\nu}+\nu_0)q_{\rho}+q_{\nu} (3d_{\rho}+\rho_0)+b (3d_{\rho}\nu_0+3d_{\nu}\rho_0+2\nu_0\rho_0)+3 b^2d_{\rho}d_{\nu} s\right)\\
& \displaystyle  \left. -d_{\nu} e^{b s} (q_{\rho}+3 b\rho_0+3 bq_{\rho} s)-q_{\rho} (\nu_0+2q_{\nu} s)-b (\nu_0+q_{\nu} s) \big(3d_{\rho} e^{a s}+2 (\rho_0+q_{\rho} s)\big)\right]\\
& \displaystyle   +a^2 b^2 \Big[-2q_{\rho}q_{\nu}+e^{(a+b) s} \Big(2q_{\rho}q_{\nu}+b^2 \big(3d_{\nu}\rho_0+\nu_0 (d_{\rho}+\rho_0)\big)+b \big((3d_{\nu}+\nu_0)q_{\rho} \\
& \displaystyle   +q_{\nu} (d_{\rho}+\rho_0)\big)+b^3d_{\rho}d_{\nu} s\Big)-3 bd_{\nu} e^{b s} (q_{\rho}+b\rho_0+bq_{\rho} s)-b^2 \left(d_{\rho} e^{a s}+\rho_0+q_{\rho} s\right) (\nu_0+q_{\nu} s)\\
& \displaystyle  -b \big(\nu_0q_{\rho}+q_{\nu} \left(d_{\rho} e^{a s}+\rho_0+2q_{\rho} s\right)\big)\Big]+
       a^4 b \left[e^{(a+b) s} \big(3 bd_{\rho}\nu_0+3d_{\rho}q_{\nu}+b (d_{\nu}+\nu_0)\rho_0  \right. \\
& \displaystyle \left. +3 b^2d_{\rho}d_{\nu} s\big)-b (\rho_0+q_{\rho} s) \left(d_{\nu} e^{b s}+\nu_0+q_{\nu} s\right)-
           3d_{\rho} e^{a s} \big(q_{\nu}+b (\nu_0+q_{\nu} s)\big)\right]\bigg\}\Bigg\}^2ds dt.\\
\end{array}
}}
$$

Note that for this general case, it is not possible to obtain an explicit expression for
$\eta_2^2(T)$, so that it can be approximated using an appropriate numerical quadrature formula.

\section{Market data}\label{Append:MarketData}

\begin{table}[h!]
\begin{footnotesize}
	\begin{center}
		\begin{tabular}{c c c c}
		\textbf{Time} & \textbf{Year fraction, $T$} & \textbf{Interest rate, $r$} & \textbf{Dividend yield, $y$} \\
		\hline\hline
		$3$ months & $0.2438$ & $1.4198$ \% & $1.5620$ \% \\
		\hline
		$6$ months & $0.4959$ & $1.2413$ \% & $2.9769$ \% \\
		\hline
		$12$ months & $1$ & $1.0832$ \% & $1.9317$ \% \\
		\hline
		$24$ months & $2$ & $1.0394$ \% & $1.8610$ \% \\
		\hline
		\end{tabular}
		\caption{EURO STOXX 50 (Dec. 2011). Spot value $S_0 = 2311.1$ \euro. Interest rates and dividend yields.}
		\label{eurostoxx1}
	\end{center}
\end{footnotesize}
\end{table}

\begin{table}[h!]
\begin{footnotesize}
	\begin{center}
		\begin{tabular}{c c c c c}
		\textbf{$K$ (\% of $S_0$)} & \textbf{$3$ months} &
                       \textbf{$6$ months} & \textbf{$12$ months} & \textbf{$24$ months} \\
		\hline\hline

		$80$\% & $33.90$\% & $33.81$\% & $31.38$\% & $29.25$\% \\
		$82$\% & $33.47$\% & $33.37$\% & $31.06$\% & $28.98$\% \\
		 $84$\% & $33.05$\% & $32.93$\% & $30.75$\% & $28.70$\% \\
		$86$\% & $32.62$\% & $32.49$\% & $30.43$\% & $28.43$\% \\
		$88$\% & $32.21$\% & $32.07$\% & $30.11$\% & $28.16$\% \\
		$90$\% & $31.79$\% & $31.64$\% & $29.79$\% & $27.89$\% \\
		$92$\% & $31.38$\% & $31.23$\% & $29.48$\% & $27.62$\% \\
		$94$\% & $30.98$\% & $30.82$\% & $29.16$\% & $27.34$\% \\
		$96$\% & $30.58$\% & $30.42$\% & $28.85$\% & $27.08$\% \\
		$98$\% & $30.18$\% & $30.02$\% & $28.53$\% & $26.81$\% \\
		$100$\% & $29.79$\% & $29.63$\% & $28.22$\% & $26.54$\% \\
		$102$\% & $29.40$\% & $29.24$\% & $27.90$\% & $26.27$\% \\
		$104$\% & $29.01$\% & $28.87$\% & $27.59$\% & $26.00$\% \\
		$106$\% & $28.63$\% & $28.49$\% & $27.28$\% & $25.74$\% \\
		$108$\% & $28.26$\% & $28.13$\% & $26.96$\% & $25.47$\% \\
		$110$\% & $27.89$\% & $27.77$\% & $26.65$\% & $25.21$\% \\
		$112$\% & $27.52$\% & $27.42$\% & $26.34$\% & $24.95$\% \\
		$114$\% & $27.16$\% & $27.07$\% & $26.03$\% & $24.68$\% \\
		$116$\% & $26.80$\% & $26.73$\% & $25.72$\% & $24.42$\% \\
		$118$\% & $26.45$\% & $26.40$\% & $25.41$\% & $24.16$\% \\
		$120$\% & $26.10$\% & $26.07$\% & $25.10$\% & $23.90$\% \\
		\end{tabular}
		\caption{EURO STOXX 50 (Dec. 2011). Implied volatilities for each maturity with different strikes $K$ (\% of the spot $S_0$).}
		\label{eurostoxx1_vola}
	\end{center}
\end{footnotesize}
\end{table}

\begin{table}[h!]
\begin{footnotesize}
	\begin{center}
		\begin{tabular}{c c c c}
		\textbf{Time} & \textbf{Year fraction, $T$} & \textbf{Interest rate, $r$} & \textbf{Dividend yield, $y$} \\
		\hline\hline
		$3$ months & $0.2528$ & $1.3696$ \% & $0.5894$ \% \\
		\hline
		$6$ months & $0.5083$ & $1.2110$ \% & $0.6185$ \% \\
		\hline
		$12$ months & $1$ & $1.0832$ \% & $0.6907$ \% \\
		\hline
		$24$ months & $2$ & $1.0394$ \% & $0.7438$ \% \\
		\hline
		\end{tabular}
		\caption{EUR/USD (Dec. 2011). Spot value $S_0 = 1.2939$ US dollars. Interest rates and dividend yields.}
		\label{eurusd}
	\end{center}
\end{footnotesize}
\end{table}

\begin{table}[h!]
\begin{footnotesize}
	\begin{center}
		\begin{tabular}{c c|c c|c c|c c}
		\multicolumn{2}{c}{\textbf{$3$ months}} & \multicolumn{2}{c}{\textbf{$6$ months}} & \multicolumn{2}{c}{\textbf{$12$ months}} & \multicolumn{2}{c}{\textbf{$24$ months}} \\
		\textbf{$K$} & $\sigma_{market}$ & \textbf{$K$} & $\sigma_{market}$ & \textbf{$K$} & $\sigma_{market}$ & \textbf{$K$} & $\sigma_{market}$ \\
		\hline\hline
		$1.1075$ & $19.27$\% & $1.0241$ & $20.75$\% & $0.9217$ & $21.80$\% & $0.8245$ & $21.01$\% \\
		$1.1516$ & $18.58$\% & $1.0877$ & $19.90$\% & $1.0084$ & $20.82$\% & $0.9302$ & $20.11$\% \\
		$1.1817$ & $17.93$\% & $1.1316$ & $19.11$\% & $1.0693$ & $19.91$\% & $1.0063$ & $19.27$\% \\
		$1.2075$ & $16.85$\% & $1.1700$ & $17.71$\% & $1.1240$ & $18.22$\% & $1.0746$ & $17.75$\% \\
		$1.2262$ & $16.40$\% & $1.1972$ & $17.19$\% & $1.1621$ & $17.64$\% & $1.1244$ & $17.21$\% \\
		$1.2425$ & $16.00$\% & $1.2210$ & $16.70$\% & $1.1956$ & $17.09$\% & $1.1686$ & $16.70$\% \\
		$1.2570$ & $15.62$\% & $1.2423$ & $16.26$\% & $1.2257$ & $16.60$\% & $1.2089$ & $16.24$\% \\
		$1.2704$ & $15.28$\% & $1.2618$ & $15.85$\% & $1.2534$ & $16.15$\% & $1.2463$ & $15.82$\% \\
		$1.2830$ & $14.97$\% & $1.2801$ & $15.49$\% & $1.2794$ & $15.74$\% & $1.2818$ & $15.44$\% \\
		$1.2950$ & $14.70$\% & $1.2975$ & $15.17$\% & $1.3043$ & $15.39$\% & $1.3161$ & $15.11$\% \\
		$1.3066$ & $14.46$\% & $1.3145$ & $14.88$\% & $1.3286$ & $15.07$\% & $1.3500$ & $14.81$\% \\
		$1.3183$ & $14.25$\% & $1.3315$ & $14.64$\% & $1.3530$ & $14.80$\% & $1.3843$ & $14.55$\% \\
		$1.3302$ & $14.07$\% & $1.3489$ & $14.44$\% & $1.3781$ & $14.58$\% & $1.4199$ & $14.34$\% \\
		$1.3427$ & $13.93$\% & $1.3673$ & $14.28$\% & $1.4047$ & $14.41$\% & $1.4579$ & $14.16$\% \\
		$1.3563$ & $13.83$\% & $1.3872$ & $14.15$\% & $1.4339$ & $14.28$\% & $1.5000$ & $14.03$\% \\
		$1.3715$ & $13.75$\% & $1.4099$ & $14.07$\% & $1.4673$ & $14.19$\% & $1.5485$ & $13.94$\% \\
		$1.3899$ & $13.74$\% & $1.4370$ & $14.02$\% & $1.5080$ & $14.16$\% & $1.6074$ & $13.88$\% \\
		$1.4140$ & $13.80$\% & $1.4733$ & $14.09$\% & $1.5635$ & $14.26$\% & $1.6890$ & $13.94$\% \\
		$1.4510$ & $13.89$\% & $1.5300$ & $14.21$\% & $1.6514$ & $14.42$\% & $1.8203$ & $14.06$\% \\
		\end{tabular}
		\caption{EUR/USD (Dec. 2011). Implied volatilities for each maturity with different strikes $K$.}
		\label{eurusd_vola}
	\end{center}
\end{footnotesize}
\end{table}

\end{document}